\pdfoutput=1
\documentclass{article}

\usepackage{mathtools}

\usepackage{graphicx}
\usepackage{amsfonts}
\usepackage{amssymb}
\usepackage{amsmath}
\usepackage{amsthm}
\usepackage{url}
\usepackage{color}
\usepackage{multirow}
\usepackage{enumerate}
\usepackage{hyperref}
\usepackage{epstopdf}
\usepackage{tikz}
\usepackage{here}
\usepackage{bigstrut}
\usepackage{pgfplots}
\usepackage{subcaption}
\usepackage{algorithm}
\usepackage{algpseudocode}
\usepackage{adjustbox}
\usepackage{xcolor, colortbl}
\usepackage[indent]{parskip}
\usepackage{afterpage}
\usepackage{dsfont}

\DeclareMathOperator{\diag}{diag}
\DeclareMathOperator{\rank}{rank}
\DeclareMathOperator{\conv}{conv}

\DeclareMathOperator{\STAB}{STAB}
\DeclareMathOperator{\COL}{COL}
\DeclareMathOperator{\T}{TH}

\newtheorem{thm}{Theorem}
\newtheorem{lem}[thm]{Lemma}
\newtheorem{cor}[thm]{Corollary}

\newtheorem{example}[thm]{Example}

\begin{document}

\title{Practical Experience with Stable Set and Coloring Relaxations 
	}
\author{
		{Dunja Pucher} \thanks{Department of Mathematics,
			University of Klagenfurt, Austria, {\tt
				dunja.pucher@aau.at}}
		\and  {Franz Rendl}\thanks{Department of Mathematics,
			University of Klagenfurt, Austria,
			{\tt franz.rendl@aau.at} }}
	
	%\date{\today}
	\date{}
	
	\maketitle
	\begin{abstract}
		
	The stable set problem and the graph coloring problem are classes of NP-hard optimization problems on graphs. It is well known that even near-optimal solutions for these problems are difficult to find in polynomial time. The Lovász theta function, introduced by Lovász in the late 1970s, provides a powerful tool in the study of these problems. It can be expressed as the optimal value of a semidefinite program and serves as a relaxation for both problems. Over the years, considerable effort has been devoted to investigating additional cutting planes to strengthen these relaxations. In our work, we use these models and consider classes of cutting planes based on cliques, odd cycles, and odd antiholes contained in the underlying graph. We demonstrate that identifying such violated constraints can be done efficiently and that they often lead to significant improvements over previous bounds. 
    \\

		\noindent Keywords: Lovász Theta Function, Stability Number, Chromatic Number, Semidefinite Programming
		
	\end{abstract}
	
\section{Introduction}
	
We consider two fundamental combinatorial optimization problems: the stable set problem and the graph coloring problem. 
Given a graph $G$, a stable set is a subset of vertices such that no two vertices in the subset are adjacent.
The cardinality of a maximum stable set is called the 
stability number of $G$ and is denoted by $\alpha(G)$. The
complement of $G$ is denoted by $\overline{G}$. Thus, stable sets in
$G$ are in one-to-one correspondence with complete subgraphs, also
called cliques, in $\overline{G}$. The number
$\omega(\overline{G}) = \alpha(G)$ denotes the cardinality of a
maximum clique in $\overline{G}$. 
        
A graph coloring is an assignment of colors to the 
vertices of a graph $G$ such that adjacent vertices
receive distinct colors. 
This can be seen as a partitioning the vertices of the graph $G$ into
$k$ stable sets $S_{1}, \ldots, S_{k}$, since such a partition
provides a $k$-coloring of $G$ by assigning the color $i$ to all
vertices in $S_{i}$ for $i=1, \ldots, k$. 
The smallest number $k$ such that $G$ has a $k$-partition into
stable sets is called the chromatic number of $G$ and is denoted by
$\chi(G)$. Partitioning of the vertex set $V$ into complete subgraphs is called clique partitioning, and the clique partition number is the smallest $k$ such that the graph has a clique partition into $k$ sets containing all edges.
Therefore, the clique partition number of $\overline{G}$  coincides
with the chromatic number of $G$. Thus, it is legitimate to consider
only the stable set and the coloring problem, as done here.  

The stable set problem and the graph coloring problem are both contained in Karp's original list of NP-complete problems~\cite{Karp72} from 1972.
Hence, unless P~=~NP, no polynomial time algorithm can solve these
problems exactly.

By definition, the inequality $\alpha(G) \leq \chi(\overline{G})$ holds for these NP-hard graph parameters. A well-known semidefinite relaxation for these parameters was introduced by Lovász in \cite{Lov:79}. The Lovász theta function, which can be computed to a fixed precision in polynomial time (see, for instance, Vandenberghe and Boyd \cite{VanBoyd}), provides a bound that separates the stability number of a graph $G$ and the chromatic number of its complement $\overline{G}$
\begin{align}\label{def-theta}
\alpha(G) \leq \vartheta(G) \leq \chi(\overline{G}).
\end{align}
	
Since its introduction in 1979, the Lovász theta function has been
extensively studied. Several authors considered various strengthening that yield tighter bounds on $\alpha(G)$ and $\chi(G)$. A standard way to strengthen a relaxation
is to add cutting planes; see, for instance, Fortet \cite {For:60}
and Barahona, J\"unger, and Reinelt \cite{BJR:89}.
This approach has been exploited, for instance, in strengthenings proposed by Schrijver~\cite{Schrijver}, Lovász and Schrijver~\cite{LovSch}, Szegedy~\cite{Szegedy}, and Meurdesoif~\cite{Meurdesoif}. 

Another way to strengthen the Lovász theta function is to use a hierarchical approach. This systematic procedure involves strengthening a relaxation by adding additional variables and constraints, and it usually unfolds in multiple levels, each providing a tighter relaxation than the previous one. The most prominent hierarchies based on semidefinite programming (SDP) are the ones from Sherali and Adams~\cite{Sherali1990}, Lovász and Schrijver~\cite{LovSch}, and Lasserre~\cite{Lasserre}. Laurent investigated the computational application of these hierarchies to the stable set polytope in~\cite{Laurent2003}.

However, any strengthening of the Lovász theta function inevitably increases computational complexity. To address this challenge, our research focuses on identifying subgraphs and structures where the separation problem can be efficiently solved and where the addition of valid inequalities into the semidefinite program (SDP) for computing the Lovász theta function can be executed within relatively short running times. In particular, in Pucher and Rendl~\cite{PuchRen2023}, we derived several inequalities valid for specific subgraphs containing maximal cliques. Our computational results indicated that incorporating these inequalities into the SDP formulation for computing the Lovász theta function significantly improves the bounds on both $\alpha(G)$ and $\chi(G)$. 

In this work, we generalize some of the results presented in~\cite{PuchRen2023} and focus on subgraphs with stability number two. Moreover, we also investigate subgraphs induced either partially or entirely by odd cycles and odd antiholes, and present several inequalities valid for the stable set and graph coloring problems. Finally, we conduct extensive computational experiments to demonstrate the effectiveness of these new inequalities, alongside those we introduced in~\cite{PuchRen2023}, in comparison to the previous state of the art.
	
The rest of this paper is organized as follows. In Section~\ref{preliminaries}, we review the literature relevant for our work. Then, in Sections~\ref{Section_inequalities} and ~\ref{inequalities_coloring}, we present new inequalities valid for the stable set and graph coloring problems, respectively. Section~\ref{Computational_results} is devoted to computational experiments. Finally, we conclude with a short discussion of the results in Section~\ref{Conclusion}.
	
We close this section with some words on notation and the terminology used. The vector of all-ones of length $n$ is denoted by $e$, and we write $J = e e^\top$ for the matrix of all-ones. Let $0_n$ be the zero matrix and $I_n$ the identity matrix of order~$n$. We denote by $e_i$ the column $i$ of the matrix $I_n$. Furthermore, we set
\begin{align*}
	E_i &= e_i e_i^T \\
	E_{ij} &= (e_i + e_j)(e_i + e_j)^T. 
\end{align*}

Throughout this paper, $G = (V(G), E(G))$ denotes a simple undirected graph with $\vert V(G) \vert = n$ vertices and $\vert E(G) \vert = m$ edges. If the graph is clear from the context, we also write $G = (V, E)$. A graph $G'$ is a subgraph of $G$
if $V(G') \subseteq V(G)$ and $E(G') \subseteq E(G)$. 
Let $V' \subseteq V(G)$. We say that the subgraph $G'$ of $G$ is induced by $V'$ if $E(G')$ consists of all edges
$\{i,j\}$ with $i,j \in V'$. The set of vertices adjacent to a vertex~$i$ is called the open neighborhood of~$i$ in~$G$ and is denoted by~$N(i)$. The closed neighborhood of a vertex~$i$, denoted~$N[i]$, contains its open neighborhood and the vertex~$i$ itself. Thus,~$N[i] = N(i) \cup \{i\}$.

A cycle in a graph is a sequence of distinct vertices~$i_1, i_2, \dots, i_k$ with~$k \geq 3$ such that each consecutive pair~$\{i_i, i_{i+1}\}$ is an edge in the graph for~$i = 1, \dots, k-1$, and the last vertex~$i_k$ is adjacent to the first vertex~$i_1$. The length of a cycle is the number of its vertices. A chord in a cycle is an edge that connects two non-consecutive vertices~$i$ and~$j$ of the cycle, where~$i \neq j$ and~$i$ and~$j$ are not adjacent in the cycle sequence. A hole in a graph is a chordless cycle of length at least four. An antihole of length~$k$ is the complement of a hole of length~$k$. An odd cycle is a cycle with an odd number of vertices. An odd hole is an odd cycle with no chords, and an odd antihole is the complement of an odd hole.

\section{Preliminaries}\label{preliminaries}

This section outlines the theoretical preliminaries necessary for this work. In Section~\ref{Formulations}, we formulate the stable set and graph coloring problems as integer programs and demonstrate how to obtain semidefinite programming relaxations for these problems. Specifically, we show that the Lovász theta function provides both upper bounds on $\alpha(G)$ and lower bounds on $\chi(\overline{G})$. Section~\ref{strengthening} then reviews existing enhancements of the Lovász theta function towards $\alpha(G)$ and $\chi(G)$, and outlines the primary objectives of this work in that context.

\subsection{Stable Sets, Colorings, and the Lovász  Theta Function}\label{Formulations}
Both the stable set problem and the graph coloring problem can be formulated as integer optimization problems. The standard formulation for the stable set problem goes back to Padberg~\cite{Padberg1973OnTF} and is given as follows. Let $x_i$ be a binary variable indicating whether the vertex $i$ is contained in a stable set ($x_i = 1$) or not ($x_i = 0$). Then, the stability number of a graph $G$ is the optimal value of the integer linear program
\begin{align}\label{stable_set_IP}
  \alpha(G) ~= ~ \max \{ e^Tx \colon x_i + x_j \leq 1 ~\forall
  \{i, j\} \in E(G), ~x_i \in\{0, 1\} ~\forall i \in V(G)\}.
\end{align}

An alternative formulation is obtained by introducing quadratic constraints; see, for instance, Ben-Tal and Nemirovski~\cite {Ben-Tal2001}. Let $x_i$ be a binary variable as in formulation~\eqref{stable_set_IP}. Then, the binary constraint can be replaced by $x_i ^2 = x_i$, and the edge inequalities $x_i + x_j \leq 1$ can be written as $x_ix_j = 0$. Thus, 
\begin{align}\label{stable_set_IP_2}
	\alpha(G) ~=~ \max \{ e^Tx \colon x_ix_j = 0 ~\forall \{i, j\} \in E(G), ~x_i^2 = x_i ~\forall i \in V(G)\}.
\end{align}

Since the computation of $\alpha(G)$ is NP-hard, one usually considers a suitable relaxation. A simple linear programming relaxation of the
formulation~(\ref{stable_set_IP}) can be obtained by omitting the
binary constraint and allowing
\begin{align}\label{nonneg}
  0 \leq x_i \leq 1 ~~\forall i \in V(G).
\end{align}

However, as Nemhauser and Trotter noted in~\cite{Nemhauser1975}, this linear relaxation is relatively weak. To strengthen it, one can add inequalities that are valid for the convex hull of the set of incidence vectors of stable sets in $G$. We denote the set of these vectors as
\begin{align*}
\mathcal{S}(G) = \{s \in \{0,1\}^{n}: s
\mbox{ incidence vector of some stable set in }G \}.
\end{align*}
The stable set polytope
$\STAB(G)$ is then defined as
\begin{align*}
\STAB(G) = \conv \{s \colon s \in \mathcal{S}(G)\}.
\end{align*}

The study of valid inequalities for $\STAB(G)$ was initiated by Padberg in~\cite{Padberg1973OnTF}. For instance, if $C \subseteq V(G)$ denotes the vertex set of an odd cycle in $G$, then clearly at
most $\frac{1}{2}(|C|-1)$ vertices of $C$ can be contained in any stable set of $G$. Moreover, if $A \subseteq V(G)$ denotes the vertex set of inducing an odd antihole in $G$, then at most two of its vertices may be contained in any stable set. Similarly, having a clique $Q \subseteq V(G)$ in $G$, at most one of its vertices may be contained in any stable set. Thus, as shown in~\cite{Padberg1973OnTF}, the following classes of inequalities are valid for $\STAB(G)$
\begin{alignat}{4}
\textrm{odd cycle constraints:} \quad \quad & \sum_{i \in C} x_i &&
\leq \frac{1}{2}(\vert C \vert - 1) && \quad \textrm{ for all odd
  cycles $C$ in }G \label{odd_cycles_in} \\
  \textrm{odd antihole constraints:} \quad \quad & \sum_{i \in  A} x_i &&
\leq 2 && \quad \textrm{ for all odd
  antiholes $A$ in }G \label{odd_antihole_in} \\
\textrm{clique constraints:} \quad \quad & \sum_{i \in Q} x_i &&
\leq 1 && \quad \textrm{ for all cliques $Q$ in } G \label{clique_standard}.
\end{alignat}

In~\cite{GLS_1986}, Grötschel, Lovász, and Schrijver demonstrated that any linear function subject to the odd cycle constraints~\eqref{odd_cycles_in} can be optimized in polynomial time. However, optimizing a linear function with the clique constraints~\eqref{clique_standard} is NP-hard; see, for instance, Nemhauser and Wolsey~\cite{Nemhauser1988}.

Besides the mentioned inequalities, many other inequalities valid for $\STAB(G)$ are known. For the sake of brevity, we list only some of these classes: odd wheel inequalities~\cite{Padberg1973OnTF}, web and antiweb inequalities introduced by Trotter in~\cite{Trotter}, rank inequalities introduced by Nemhauser and Trotter in~\cite{Nemhauser1974}; some further inequalities can be found in works from Cheng and Cunningham~\cite{Cheng1997}, C{\'a}novas, Landete, and Mar{\'{\i}}n~\cite{Canovas2000}, and Cheng and de Vries~\cite{Cheng2002}.

Another way to obtain upper bounds on $\alpha(G)$ is to consider the Lovász theta function $\vartheta(G)$, introduced by Lovász~\cite{Lov:79} in his seminal paper from 1979. A semidefinite program to compute the Lovász theta function can be formulated in several ways, as demonstrated, for instance, by Grötschel, Lovász, and Schrijver in~\cite{GLS_1988} and by Knuth in~\cite{Knuth1994}. Here, we state the formulation that was given by Lovász and Schrijver in~\cite{LovSch} and that can be directly derived by considering the integer programming formulation with quadratic constraints~\eqref{stable_set_IP_2}, as outlined, for instance, in Rendl~\cite{Rendl2010}.

Suppose that a vector $x \in \{0, 1\}^n$ is feasible for the formulation~\eqref{stable_set_IP_2} and consider the matrix $X = xx^T \in \{0, 1\}^{n \times n}$. By definition, $X$ is positive semidefinite, and since $x$ is a binary vector, $x_i^2 = x_i$ holds for every $i \in V(G)$. Hence, the main diagonal of $X$ is equal to $x$. Furthermore, for every edge $\{i, j\} \in E$, the entry $X_{ij}$ equals zero. Nevertheless, the constraint $X = xx^T$ is not convex. Relaxing $X - xx^T = 0$ to $X - xx^T \succeq 0$ and using the Schur complement
\begin{align*}
X - xx^T \succeq 0 ~\Leftrightarrow~ \begin{pmatrix} X & x \\ x^T & 1 \end{pmatrix} \succeq 0,
\end{align*}
we obtain the following SDP to compute the Lovász theta function
\begin{align}\label{theta_stable_set}
  \vartheta(G) ~=~ \max \Biggl\{e^T x \colon
\begin{pmatrix} X & x \\ x^T & 1
\end{pmatrix} \succeq 0, ~ \mathrm{diag}(X) = x,
 ~ X_{ij} = 0 ~\forall \{i, j\} \in E(G)\Biggr\}.
\end{align}

In~\cite{GLS_1988}, Grötschel, Lovász, and Schrijver introduced the theta body of a graph $G$ as the convex body $\T(G) \subseteq \mathbb{R}^n$ given by the feasible region of the SDP to compute the Lovász theta function. Hence, considering the formulation~\eqref{theta_stable_set}, we obtain
\begin{align*}
  \T(G) = \Biggl\{x \in \mathbb{R}^n : \exists X \in \mathbb{S}^n :
\begin{pmatrix} X & x \\ x^T & 1
\end{pmatrix} \succeq 0, ~ \mathrm{diag}(X) = x,~ X_{ij} = 0 ~\forall \{i, j\} \in E(G)\Biggr\}.
\end{align*}

The theta body has the advantage that we can optimize a linear function over $\T(G)$ in polynomial time with an arbitrary fixed precision using semidefinite programming. Furthermore, by definition, we can see that the theta body contains the stable set polytope of $G$, i.e.,
\begin{align*}
    \STAB(G) \subseteq \T(G) \subseteq \mathbb{R}^n.
\end{align*}
Additionally, feasibility for the problem in \eqref{theta_stable_set} 
also implies that $0 \leq x_i \leq 1$. Moreover, Grötschel, Lovász, and Schrijver~\cite{GLS_1988} showed that $\T(G)$ satisfies all clique constraints. Therefore, $\vartheta(G)$ often provides a tight upper bound on the stability number $\alpha(G)$. This bound is generally stronger than those obtained from linear relaxations, as shown, for instance, by Balas, Ceria, Cornu{\'e}jols, and Pataki in~\cite{Balas1996}.

In the next section, we will discuss the ways to strengthen the Lovász theta function $\vartheta(G)$ as an upper bound on $\alpha(G)$ by adding valid inequalities into the SDP for its computation. Clearly, inequalities valid for $\STAB(G)$ can be added into~\eqref{theta_stable_set} in order to obtain better bounds. Nevertheless, we will also consider some inequalities valid for the matrix sets associated with $\STAB(G)$. Therefore, it will turn out to be useful to work with the matrix version of $\mathcal{S}(G)$. We denote it by $\mathcal{S}^{2}(G)$ and define it as
\begin{align*}
\mathcal{S}^{2}(G) =  \{ ss^T : s \in \mathcal{S}(G) \}.
\end{align*}
Its convex hull is
\begin{align*}
\STAB^{2}(G) = \conv \{ \mathcal{S}^{2}(G)\}. 
\end{align*}

We now consider the graph coloring problem and show that the Lovász theta function also yields a lower bound on $\chi(\overline{G})$. To this end, we follow the presentation by Dukanovic and Rendl from~\cite{Dukanovic2007}. Recall that a $k$-coloring of a graph G on $n$ vertices is a $k$-partition $(V_1,\ldots, V_k)$ of the vertex set $V(G)$ such that each $V_i$ is a stable set in $G$. We encode the $k$-partition by the characteristic vectors $s_i$ for $V_i$. Hence, $s_i \in \{0, 1\}^n$ and for a $v \in V(G)$ we have $(s_i)_v = 1$ if $v \in V_i$, and $0$ otherwise. Furthermore, the partition property implies that $s_i \neq 0$ and $\sum_{i = 1}^k s_i = e$. The matrix 
\begin{align*}
X = \sum_{i = 1}^k s_i s_i^T
\end{align*}
is called a coloring matrix. The convex hull of the set of all coloring matrices of $G$ is denoted by
\begin{align*}
\COL(G) = \conv\{X : X \text{ is a coloring matrix of $G$}\}.
\end{align*}

Hence, a coloring matrix $X$ is the partition matrix associated to the $k$-partition $(V_1,\ldots, V_k)$. Clearly, $\rank(X) = k$ and $X$ is a symmetric $0/1$ matrix. Partition matrices can be characterized in several ways. In particular, a symmetric matrix $X \in \{0, 1\}^{n \times n}$ is a $k$-partition matrix if and only if $ \diag(X) = e$ and $t X - J \succeq 0$ for all $t \geq k$. Consequently, the chromatic number of a graph $G$ is the optimal value of the following SDP in binary variables
\begin{align}\label{coloring_IP}
\chi(G) ~=~ \min \{t \colon t X - J \succeq 0, ~ \mathrm{diag}(X) = e,
   ~ X_{ij} = 0 ~\forall \{i, j\} \in E(G), ~X \in \{0, 1\}^{n \times n} \}.
\end{align}

Nevertheless, solving~\eqref{coloring_IP} is NP-hard. A tractable relaxation can be obtained by leaving out the binary condition. Furthermore,
\begin{align*}
t X - J \succeq 0 ~\Leftrightarrow~ \begin{pmatrix} t & e^T \\ e & X \end{pmatrix} \succeq 0
\end{align*}
by the Schur complement. Hence, we obtain the following semidefinite relaxation that yields a lower bound on the chromatic number $\chi(G)$
\begin{align}\label{theta_coloring}
\vartheta (\overline{G}) ~=~ \min \Biggl\{t \colon
\begin{pmatrix} t & e^T \\ e & X
\end{pmatrix} \succeq 0, ~ \mathrm{diag}(X) = e,
   ~ X_{ij} = 0 ~\forall \{i, j\} \in E(G)\Biggr\}.
\end{align}

Note that the optimal value of~\eqref{theta_coloring} is indeed the Lovász theta function $\vartheta(\overline{G})$. This can be demonstrated by considering the dual of the SDP problem~\eqref{theta_stable_set} applied to the complement of $G$. In particular, this proves the Lovász sandwich theorem, as stated in equation~\eqref{def-theta} in the introduction of this paper. 

\subsection{Strengthening of the Lovász Theta Function}\label{strengthening}

In the previous section, we have seen that the parameter $\vartheta(G)$ as an upper bound on $\alpha(G)$ is given as the optimal value for the semidefinite optimization problem \eqref{theta_stable_set}. There are some obvious strengthenings towards $\alpha(G)$. The first strengthening $\vartheta^\prime(G)$ of the Lovász theta function was given by Schrijver~\cite{Schrijver} by adding the non-negativity constraints
\begin{align}\label{non-n}
  X \geq 0.
\end{align}  
Gruber and Rendl~\cite{Gru:03} proposed a further strengthening
$\vartheta^{\prime \Delta}(G)$ which is obtained by including the triangle inequalities 
\begin{align}
X_{ik} + X_{jk} \leq X_{ij} + x_k, &\quad \forall \{i, j, k\} \subseteq V(G) \label{triangle_inequalities_1} \\
 x_{i}+x_{j}+x_{k} \leq 1 + X_{ik} + X_{jk} + X_{ij}, &\quad \forall \{i, j, k\} \subseteq V(G) \label{triangle_inequalities_2}
\end{align}
into the SDP \eqref{theta_stable_set}, yielding a chain of relaxations
\begin{align*}
\alpha(G) \leq \vartheta^{\prime \Delta}(G) \leq \vartheta^ \prime (G) \leq  \vartheta(G).
\end{align*}

The validity of the triangle inequalities~\eqref{triangle_inequalities_1} and~\eqref{triangle_inequalities_2} for the stable set polytope was first observed by Padberg in~\cite{PadbergBP}. In that work, Padberg introduced triangle inequalities valid for the boolean quadric polytope, which is the convex hull of all matrices $X$ of the form $xx^T$ where $x \in \{0, 1\}^n$. Nevertheless, the polytope $\STAB(G)$ is derived by taking the face of the boolean quadric polytope defined by the equations $X_{ij} = 0$ for all $\{i, j\} \in E(G)$ and then projecting it onto the non-quadratic space. Consequently, the triangle inequalities valid for the boolean quadric polytope also hold for the stable set polytope.

Additionally, it is important to mention the relaxation proposed by Lovász and Schrijver in~\cite{LovSch}, which includes both the non-negativity constraint~\eqref{non-n} and a subset of the triangle inequalities~\eqref{triangle_inequalities_1} and~\eqref{triangle_inequalities_2}. This relaxation is derived using the Lift-and-Project operator introduced in the same work. While it does not encompass all triangle inequalities, making it slightly weaker than $\vartheta^{\prime \Delta}(G)$, it does satisfy all odd cycle inequalities~\eqref{odd_cycles_in}. Therefore, incorporating the non-negativity constraint~\eqref{non-n} along with the triangle inequalities~\eqref{triangle_inequalities_1} and~\eqref{triangle_inequalities_2} into the SDP~\eqref{theta_stable_set} ensures that all odd cycle inequalities~\eqref{odd_cycles_in} are satisfied.

Several strengthenings of $\vartheta(\overline{G})$
towards $\chi(G)$ have been proposed
and are quite similar to the strengthening of
$\vartheta(G)$ towards $\alpha(G)$. 
Szegedy~\cite{Szegedy} introduced the non-negativity constraint 
\begin{align}\label{non-n-col}
X \geq 0
\end{align}leading to $\vartheta^{+}(\overline{G})$. 
Meurdesoif~\cite{Meurdesoif} introduced an even tighter
bound $\vartheta^{{+\Delta}}(\overline{G})$,
which is obtained by adding the triangle inequalities 
\begin{align}\label{triangle_inequalities_coloring}
X_{ij} + X_{jk} \leq X_{ik} + 1, \quad \forall \{i, j, k\} \subseteq V(G)
\end{align}
into the SDP \eqref{theta_coloring}, resulting again in achain of relaxations
\begin{align*}
  \vartheta(\overline{G})  \leq
  \vartheta^+(\overline{G}) \leq
  \vartheta^{+\Delta}(\overline{G}) \leq \chi(G).
\end{align*}

We have just seen several rather obvious ways to tighten
the Lovász theta function either towards $\alpha(G)$ or
towards $\chi(G)$.
In particular, any systematic tightening towards $\alpha(G)$ or $\chi(G)$ amounts to identifying a subset
of linear inequalities which provide a reasonably accurate description
of the underlying convex hull.

A general cutting plane theory along these lines is given by the class of Chv\'atal-Gomory cuts, see for instance
Caprara and Fischetti~\cite{caprara-fischetti}. However, identifying violated Chv\'atal-Gomory cuts is NP-hard in general. This motivates the search for problem-specific cutting planes that exploit the combinatorial structure of the underlying problem. Gr\"otschel and Wakabayashi~\cite{gro-wak-1} consider the clique partition problem in $G$ which corresponds
to coloring in $\overline{G}$ and describe several
classes of cutting planes and investigate their
practical impact, see \cite{gro-wak-2}. 
Strengthening the relaxation for the stable set problem by adding odd cycle and triangle inequalities was studied
by Gruber and Rendl in~\cite{Gru:03}. Computational experiments regarding strengthening both problems
by adding triangle inequalities were done by Dukanovic and Rendl in~\cite{Dukanovic2007} and Battista and De Santis in \cite{battista}. 

One of the most notable results in bounding the stability number $\alpha(G)$ and the chromatic number $\chi(G)$ based on the semidefinite programming relaxation was obtained recently by Gaar and Rendl in~\cite{Gaa:20}. These bounds were achieved by employing the exact subgraph hierarchy introduced by Anjos, Adams, Rendl, and Wiegele in~\cite{AARW:15}. This hierarchy is particularly well-suited for classes of NP-hard problems based on graphs, where the projection of the problem onto a subgraph retains the same structure as the original problem. 

In the context of the stable set and graph coloring problems, the authors proposed that for certain subsets of vertices $I \subseteq V$, the corresponding subgraph $G_I$ should be contained in $\STAB^2(G_I)$ or $\COL(G_I)$ restricted to the respective subgraph. The exact subgraph hierarchy begins by computing the Lovász theta function at the first level of the hierarchy and then progressively considers all subsets of vertices with increasing cardinalities. While this method produces tight bounds on $\alpha(G)$ and $\chi(G)$, it is associated with significant computational challenges due to the complexity involved in describing the convex hull of the underlying subgraph.

Motivated by these challenges, our research presented in~\cite{PuchRen2023} focused on identifying certain subgraphs for which the description of $\STAB^2(G)$ and $\COL(G)$ is relatively simple and can be efficiently integrated into the SDP for computing the Lovász theta function. Specifically, we examined subgraphs induced by the join of two cliques and showed the following result for the stable set problem.

\begin{lem}\label{join_of_cliques}
Let~$G = (V, E)$ be a graph, let~$Q_1, Q_2\subseteq V$ induce two cliques in~$G$ such that~$Q_1 \cap Q_2 = \emptyset$ and such that~$\alpha(G_{Q_1 \cup Q_2}) = 2$. Now let~$\tilde{Q} = Q_1 \cup Q_2$. Then,~$X_{\tilde Q} \in \STAB^2(G_{\tilde Q})$ if and only if the constraints
    \begin{alignat}{3}
    X_{ii} &\geq 0 &&  \quad \forall i \in \tilde{Q} \nonumber \\
    X_{ij} &= 0 && \quad \forall \{i, j\} \in E(G_{\tilde{Q}}) \nonumber \\
    X_{ij} &\geq 0 && \quad \forall \{i, j\} \in \overline{E}(G_{\tilde{Q}})  \nonumber \\
    \sum_{\substack{\{i, j\} \in \overline{E}(G_{\tilde{Q}})}} X_{ij} &\leq X_{ii} && \quad \forall i \in \tilde Q \label{cliques_1} \\
    \sum_{i \in \tilde{Q}} X_{ii} &\leq 1 + \sum_{\{i, j\} \in \overline{E}(G_{\tilde{Q}})} X_{ij} \label{cliques_3}
    \end{alignat}
    hold.
\end{lem} 

Now, the statement of Lemma~\ref{join_of_cliques} can be further refined. Assume that~$Q_2$ is a single vertex. Then, the next statement holds.

\begin{cor}\label{cor_clique_and_vertex}
Let~$G = (V, E)$ be a graph, let~$Q \subseteq V$ induce a clique in~$G$, and let~$k \in V \setminus Q$ be a vertex such that~$\alpha(G_{Q \cup \{k\}}) = 2$. If~$X_{Q \cup \{k\}} \in \STAB^2(G_{Q \cup \{k\}})$, then
    \begin{alignat}{3}
    \sum_{i \in Q} X_{ik} &\leq X_{kk}.  \label{cliques_final}
    \end{alignat}
\end{cor}
 
 Similarly, for the graph coloring problem, we established the following statement.

\begin{lem}\label{clique_vertex_clique}
Let~$G = (V, E)$ be a graph, let~$Q \subseteq V$ be a clique in~$G$, and let~$k \in V$ be a vertex such that~$k \notin Q$. Then,~$X_{Q \cup \{k\}} \in \COL(G_{Q \cup \{k\}})$ if and only if the constraints
    \begin{alignat}{3}
    X_{ii} &= 1 &&\quad \forall i \in Q \cup \{k\} \label{clique_vertex_one}\\
    X_{ij} & = 0 &&\quad \forall ~\{i, j\} \in E(G_{Q \cup \{k\}}) \label{clique_vertex_two}\\
    X_{ij} &\geq 0 &&\quad \forall ~\{i, j\} \in \overline{E}(G_{Q \cup \{k\}}) \label{clique_vertex_three}\\
    \sum_{i \in Q} X_{ik} &\leq 1 &&\label{clique_vertex_four}
    \end{alignat}
    are satisfied.  
\end{lem}

 However, it is important to note that Giandomenico, Letchford, Rossi, and Smriglio arrived at a similar result while studying the strengthening of the linear relaxation for the stable set problem in~\cite{Giandomenico2009}. Specifically, by considering the Lift-and-Project operator introduced by Lovász and Schrijver in~\cite{LovSch}, they derived the same class of inequalities~\eqref{cliques_3} and \eqref{cliques_final}. However, these inequalities were only applied to the linear relaxation. Further computational studies of these inequalities within the context of linear relaxation were conducted by Giandomenico, Rossi, and Smriglio in~\cite{Giandomenico2013}. To the best of our knowledge, paper~\cite{PuchRen2023} is the first to explore the application of inequalities~\eqref{cliques_3} and \eqref{cliques_final} in the SDP setting. Notably, computational results from~\cite{PuchRen2023} demonstrate that incorporating these inequalities into the SDP to compute the Lovász theta function significantly improves bounds on both $\alpha(G)$ and $\chi(G)$.

This success has motivated our current study, where we explore additional inequalities valid for $\STAB^2(G)$ and $\COL(G)$ by focusing on specific subgraphs of $G$ that are easily separable. This time, these new inequalities are valid for subgraphs induced by odd cycles and odd antiholes, either partially or entirely.

Nevertheless, before presenting these results, we revisit the findings from~\cite{PuchRen2023}, and note that subgraphs induced by the join of two maximal cliques have a stability number of two. Inspired by this observation, we extend our analysis to all subgraphs with stability number two.

\section{Valid Inequalities for the Stable Set Problem}\label{Section_inequalities}

    In this section, we propose several inequalities that are valid for the stable set problem. In Section~\ref{section_convex_hulls}, we focus on subgraphs with stability number two. Then, in Section~\ref{section_odd_cycles}, we investigate subgraphs induced by odd cycles. We begin with~$5$-cycles and extend our analysis to odd cycles of arbitrary length. A similar approach is applied in Section~\ref{section_antiholes}, where we present results for subgraphs containing odd antiholes.

\subsection{On the Convex Hull of Certain Subgraphs}\label{section_convex_hulls}
   
    As a first step, we focus on subgraphs whose convex hulls can be described in a relatively straightforward manner. In particular, the convex hull of subgraphs with stability number two admits a particularly simple description, as the following result demonstrates.
    
    \begin{thm}\label{stability_number_2}
    Let~$G = (V, E)$ be a graph and let~$I\subseteq V$ such that~$\alpha(G_I) = 2$. Then,~$X_I \in \STAB^2(G_I)$ if and only if the constraints
    \begin{alignat}{3}
    X_{ii} &\geq 0 &&\quad \forall i \in I \label{minus_one}\\
    X_{ij} & = 0 &&\quad \forall ~\{i, j\} \in E(G_I) \label{minus_two}\\
    X_{ij} &\geq 0 &&\quad \forall ~\{i, j\} \in \overline{E}(G_I) \label{one}\\
    \sum_{\{i,j \} \in \overline {E}(G_I)} X_{ij} &\leq X_{ii} &&\quad \forall ~i \in I \label{two}\\
    \sum_{i \in I} X_{ii} &\leq 1 + \sum_{ \{i, j\} \in \overline{E}(G_I)} X_{ij} \label{three}
    \end{alignat}
    are satisfied.
    \end{thm}	
    
    \begin{proof}
    Without loss of generality, let~$V = \{ 1, \ldots, n \}$ and~$I = \{1, \ldots, k\}$. We proceed by proving both directions of the equivalence separately. 
    
    First, we assume that~$X_I \in \STAB^2(G_I)$. Then, according to the definition and since $\alpha(G_I) = 2$, there exist non-negative coefficients~$\lambda$,~$\mu_i$ for all~$1 \leq i \leq k$ and~$\nu_{ij}$ for all~$1\leq i < j \leq k$ such that
	\begin{align}
		\lambda + \sum_{i = 1}^{k} \mu_i + \sum_{i = 1}^{k - 1} \sum_{j = i + 1}^{k} \nu_{ij} = 1 \label{sum_is_one}
	\end{align}
	and 
	\begin{align}
		X_I = \lambda 0_k + \sum_{i = 1}^{k} \mu_i E_i + \sum_{i = 1}^{k - 1} \sum_{j = i + 1}^{k} \nu_{ij} E_{ij} \label{matrix_X},
	\end{align}	
    where~$\nu_{ij} = 0$ for all pairs~$1 \leq i < j \leq k$ such that~$\{i, j\} \in E(G_I)$. 

    Thus, equating the coefficients in~\eqref{matrix_X} yields
    \begin{alignat}{3}
        X_{ij} &= \nu_{ij} &&  \quad 1 \leq i < j \leq k \label{X_ij_wo_edge} \\
        X_{ii} &= \mu_i + \sum_{\substack{j = 1 \\ j \neq i}}^k \nu_{ij} && \quad 1 \leq i \leq k \label{x_i}.
	\end{alignat}

    Since~$\nu_{ij} = 0$ for all pairs~$1 \leq i < j \leq k$ such that~$\{i, j\} \in E(G_I)$, and because~$X_{ij} = \nu_{ij}$ according to~\eqref{X_ij_wo_edge}, we obtain~\eqref{minus_two}. Furthermore, given that~$\nu_{ij} \geq 0$ for all other off-diagonal elements in~$X_I$, we obtain exactly~\eqref{one}. Moreover,~$\mu_i \geq 0$ for all~$1 \leq i \leq k$, so from~\eqref{x_i} it follows that all diagonal elements are non-negative, so~\eqref{minus_one} also holds. 

    Now, using~\eqref{X_ij_wo_edge}, we can write~\eqref{x_i} as
    \begin{align}
    X_{ii} = \mu_i + \sum_{\substack{j = 1 \\ j \neq i}}^k X_{ij} \quad 1 \leq i \leq k \label{x_i_again}.
    \end{align}
    
    Note that for a fixed~$i \in I$ we have
    \begin{align}\label{sum}
    \sum_{\substack{j = 1 \\ j \neq i}}^k X_{ij} = \sum_{\{i, j\} \in \overline{E}(G_I)} X_{ij},
    \end{align}
    because~$X_{ij} = 0$ for all~$\{i, j\} \in E(G_I)$. Given that all coefficients~$\mu_i$ are non-negative,~\eqref{x_i_again} and~\eqref{sum} yield~\eqref{two}.

    Summing up both sides of~\eqref{x_i_again} over all~$i \in I$ and noting that~$X_I$ is symmetric, we obtain
    \begin{align*}
    \sum_{i=1}^k \mu_i &= \sum_{i=1}^k X_{ii} - \sum_{i = 1}^k \sum_{\substack{j = 1 \\ j \neq i}}^k X_{ij} \\
    & = \sum_{i=1}^k X_{ii} - 2 \sum_{i = 1}^{k-1} \sum_{j = i + 1}^k X_{ij}.
    \end{align*}
    Thus, from~\eqref{sum_is_one} and keeping in mind that~$X_{ij} = 0$ for all~$\{i, j\} \in E(G_I)$, we find that
	\begin{alignat*}{2}
        \lambda &= 1 - \sum_{i = 1}^{k} \mu_i - \sum_{i = 1}^{k - 1} \sum_{j = i + 1}^{k} \nu_{ij} \\
        & = 1 - \sum_{i = 1}^k X_{ii} +  \sum_{i = 1}^{k-1} \sum_{j = i + 1}^k X_{ij} \\
        & = 1 - \sum_{i \in I} X_{ii} +  \sum_{\{i,j\} \in \overline{E}(G_I)} X_{ij}.
	\end{alignat*}
   Finally, since the coefficient~$\lambda$ is non-negative, we obtain~\eqref{three}. 

   Now assume that~\eqref{minus_one} -- \eqref{three} hold. We have to show that there exist non-negative coefficients~$\lambda$,~$\mu_i$ for all~$1 \leq i \leq k$ and~$\nu_{ij}$ for all~$1\leq i < j \leq k$ such that~\eqref{sum_is_one} and~\eqref{matrix_X} hold. We define~$\nu_{ij} = X_{ij}$ for all~$1 \leq i < j \leq k$,~$\mu_i = X_{ii} - \sum_{\substack{j = 1 \\ j \neq i}}^k X_{ij}$ for all~$1 \leq i \leq k$ and~$ \lambda = 1 - \sum_{i \in I} X_{ii} +  \sum_{\{i,j\} \in \overline{E}(G_I)} X_{ij}$. Then, from the construction, all coefficients are non-negative, and their sum equals one, so~\eqref{sum_is_one} holds. Furthermore, the matrix~$X_I$ can be expressed as~\eqref{matrix_X}. Thus,~$X_I \in \STAB^2(G_I)$, so the statement holds.
\end{proof}

    The statement of Theorem~\ref{stability_number_2} is important because it shows that if the stability number of a subgraph is two, then the convex hull of such a subgraph on~$k$ vertices can be completely described by adding~$k + 1$ inequalities, in addition to edge and non-negativity constraints. While there are various subgraphs with stability number two, such as those induced by any two non-adjacent vertices, our focus is on subgraphs of larger order~$k$ to fully leverage the potential of the given theorem. In particular, subgraphs induced by a join of two cliques demonstrate these desirable properties, as already shown in~\cite{PuchRen2023}. Thus, Lemma~\ref{join_of_cliques} established in~\cite{PuchRen2023} can be seen as a special case of Theorem~\ref{stability_number_2}.
    
    Interestingly, we observe from Lemma~\ref{join_of_cliques} that, although we consider a subgraph induced by a join of two cliques, we do not have to impose the standard clique constraints~\eqref{clique_standard} to obtain a solution that is contained in the convex hull. This leads to our next result.

\begin{lem}\label{lemma_clique_constraints_satisfied}
Let~$G$,~$Q_1$,~$Q_2$, and~$\tilde{Q}$ be as in Lemma~\ref{join_of_cliques}. Furthermore, let~$X_{\tilde Q} \in \STAB^2(G_{\tilde Q})$. Then, the constraints
    \begin{align*}
    \sum_{i \in Q_1} X_{ii}  \leq 1 ~\text{ and } ~ \sum_{j \in Q_2} X_{jj}  \leq 1 
    \end{align*}
    are satisfied.
\end{lem}

\begin{proof}
Given that~$X_{\tilde Q} \in \STAB^2(G_{\tilde Q})$, we know that~\eqref{cliques_1} holds. Summing up this constraint over all vertices~$i \in Q_1$, we obtain
\begin{align*}
\sum_{i \in Q_1} \sum_{\{i, j\} \in \overline{E}(G_{\tilde{Q}})} X_{ij} \leq \sum_{i \in Q_1} X_{ii}.
\end{align*}

Now, 
\begin{align*}
\sum_{i \in Q_1} \sum_{\{i, j\} \in \overline{E}(G_{\tilde{Q}})} X_{ij} = \sum_{i \in Q_1} \sum_{\substack{j \in Q_2 \\ \{i, j\} \in \overline{E}(G_{\tilde{Q}})}} X_{ij} = \sum_{\{i, j\} \in \overline{E}(G_{\tilde{Q}})} X_{ij},
\end{align*}
because there are no two vertices~$i, j \in Q_1$ such that~$\{i, j\} \in \overline{E}(G_{\tilde{Q}})$. Therefore, we deduce
\begin{align*}
\sum_{\{i, j\} \in \overline{E}(G_{\tilde{Q}})} X_{ij} \leq \sum_{i \in Q_1} X_{ii}.
\end{align*}

Since~\eqref{cliques_3} also holds, altogether, we obtain
\begin{align*}
\sum_{i \in Q_1} X_{ii} + \sum_{j \in Q_2} X_{jj} &= \sum_{i \in \tilde{Q}} X_{ii} \leq 1 + \sum_{\{i, j\} \in \overline{E}(G_{\tilde{Q}})} X_{ij} \leq 1 + \sum_{i \in Q_1} X_{ii}.
\end{align*}
Hence, the constraint~$\sum_{j \in Q_2} X_{jj} \leq 1$ is satisfied.

Next, we sum up constraint~\eqref{cliques_1} over all vertices~$j \in Q_2$, and following the same argumentation, we obtain that the constraint~$\sum_{i \in Q_1} X_{ii} \leq 1$ is also satisfied.
\end{proof}

    \subsection{Subgraphs Containing Odd Cycles}\label{section_odd_cycles}

We now examine subgraphs containing odd cycles and propose several inequalities valid for the stable set problem. We begin by considering subgraphs induced by odd cycles of length~$5$. Recall that any such subgraph has a stability number of at most two, which, depending on its exact stability number, allows us to potentially apply Theorem~\ref{stability_number_2} to fully describe its convex hull. However, we first present the following result, which holds for any subgraph with a stability number of at most two.

\begin{lem}\label{thm_stable_set_1}
	Let~$G = (V, E)$ be a graph, let~$I \subseteq V$ such that~$\alpha(G_I) \leq 2$. Then, the inequality 
	\begin{align}
		\sum_{\substack{i,j \in I \\ i < j}} X_{ij} \leq 1 \label{stable_set_inequality}
	\end{align}
    is valid for~$X_I \in \STAB^{2}(G_I)$.
\end{lem}

\begin{proof}
    Let~$s = (s_1, \ldots, s_n)^\top$ be the incidence vector of a stable set in~$G_I$ and define~$X_I = ss^\top$. Then, we have
    \begin{align*}
		\sum_{\substack{i,j \in I \\ i < j}} X_{ij} = \sum_{\substack{i,j \in I \\ i < j}} s_i s_j.
	\end{align*}
    
    Given that~$\alpha(G_I) \leq 2$, the largest stable set in~$G_I$ contains at most two vertices. Therefore, for all~$i, j \in I$ with~$i < j$, there is at most one product~$s_is_j$ equal one. Consequently, we have
    \begin{align*}
		\sum_{\substack{i,j \in I \\ i < j}} s_i s_j \leq 1,
	\end{align*}
    so~\eqref{stable_set_inequality} holds.
    
    Furthermore, since this holds for all incidence vectors of stable sets in~$G_I$, the inequality~\eqref{stable_set_inequality} is valid for any~$X_I \in \mathcal{S}^2(G_I)$, i.e., for any stable set matrix $ss^\top$ with $s \in \mathcal{S}(G_I)$, and thus for any~$X_I \in \STAB^{2}(G_I)$.
\end{proof}

In the context of~$5$-cycles, this leads to the following statement.

\begin{cor}\label{lemma_for_thm_stable_set_1}
	Let~$G = (V, E)$ be a graph and let~$C \subseteq V$ induce a~$5$-cycle in~$G$. Then, the inequality 
	\begin{align}
		\sum_{\substack{i,j \in C \\ i < j}} X_{ij} \leq 1 \label{stable_set_1}
	\end{align}
    is valid for~$X_C \in \STAB^{2}(G_C)$.
\end{cor}

As discussed in Section~\ref{Formulations}, another inequality valid for the~$5$-cycles is the standard odd-cycle constraint~\eqref{odd_cycles_in}. We will now show an interesting result regarding the relationship between inequalities~\eqref{stable_set_1} and~\eqref{odd_cycles_in}. More specifically, we will demonstrate that if we strengthen the Lovász theta function with~\eqref{stable_set_1}, then the standard 5-cycle inequalities~\eqref{odd_cycles_in} will be satisfied. However, the converse is not true: adding~\eqref{odd_cycles_in} into the SDP for the computation of the Lovász theta function does not imply~\eqref{stable_set_1}. Therefore, including constraints~\eqref{stable_set_1} into the formulation~\eqref{theta_stable_set} may yield tighter upper bounds on~$\alpha(G)$ than adding inequalities~\eqref{odd_cycles_in}.

\begin{lem}\label{implies}
    Let~$G = (V, E)$ be a graph, let~$C \subseteq V$ induce a~$5$-cycle in~$G$, and let~$(x, X)$ be a feasible solution for~\eqref{theta_stable_set} strengthened with the inequality~\eqref{stable_set_1}. Then,
	\begin{align*}
		\sum_{i \in C} x_i \leq 2.
	\end{align*}
\end{lem}
	
\begin{proof}
	Let~$X_C$ be the submatrix of~$X$ induced by the cycle~$C$. Since~$(x, X)$ is a feasible solution for~\eqref{theta_stable_set} strengthened with the additional inequality, the matrix~$Y = \begin{pmatrix} X_C & x_C \\ x_C^\top & 1	\end{pmatrix}$, where~$\diag(X_C) = x_C$, is positive semidefinite. Hence, for any vector~$a \in \mathbb{R}^6$, we have~$a^\top Y a \geq 0$. Now, consider the vector~$a = (1, 1, 1, 1, 1, -2)^\top$. Then,~$a^\top Y a \geq 0$ expands to
	\begin{align*}
		4 - 3 \sum_{i \in C} x_i + 2 \sum_{\substack{i,j \in C \\ i < j}} X_{ij} \geq 0,
	\end{align*}
	because~$X_C$ is symmetric. Furthermore, since $\sum_{\substack{i,j \in C \\ i < j}} X_{ij} \leq 1$ according to the assumption, we deduce
	\begin{align*}
		3 \sum_{i \in C} x_i \leq 6.
	\end{align*}
	Thus, the inequality~$\sum_{i \in C} x_i \leq 2$ holds.
\end{proof}
 
To demonstrate that adding the standard~$5$-cycle inequality~$\sum_{i \in C} x_i \leq 2$ into the SDP to compute the Lovász theta function does not imply~\eqref{stable_set_1}, we provide the following counterexample. 

\begin{example}\label{first_example}
Let~$G = (C, E(C))$ be a~$5$-cycle with the vertex set~$C = \{1, \dots, 5\}$ and the edge set $E(C) = \{\{1,2\}, \{2,3\}, \{3,4\}, \{4,5\}, \{1,5\}\}$. We define the variables as follows 
    \begin{alignat*}{3}
		X_{ii} &= 0.4 &&\quad \forall i \in C \\
        x_i &= X_{ii} &&\quad \forall i \in C \\
		X_{ij} &= 0 &&\quad \forall \{i, j\} \in E(C) \\
		X_{ij} &= 0.209 &&\quad \forall \{i, j\} \in \overline{E}(C).
    \end{alignat*}
Then,~$(x, X)$ is feasible for~\eqref{theta_stable_set}, since~$\begin{pmatrix} X & x \\ x^\top & 1 \end{pmatrix} \succeq 0$,~$X_{ij} = 0$ for all~$\{i, j\} \in E(C)$ and~$\diag(X) = x$. Furthermore, by construction,~$\sum_{i \in C} x_i = 2$. However, we find that~$\sum_{\substack{i,j \in C \\ i < j}} X_{ij} = 1.045 > 1$. 
\end{example}

As mentioned in Section~\ref{strengthening}, the triangle inequalities~\eqref{triangle_inequalities_1} and \eqref{triangle_inequalities_2} also imply the odd cycle constraints~\eqref{odd_cycles_in}. Moreover, it turns out that adding triangle inequalities into the SDP to compute the Lovász theta function implies~\eqref {stable_set_1} as well, as the next statement shows.

\begin{lem}\label{triangles_imply_sum_less_equal_one}
Let~$G = (V, E)$ be a graph, let~$C \subseteq V$ induce a~$5$-cycle in~$G$, and let $(x, X)$ be a feasible solution for~\eqref{theta_stable_set} strengthened with triangle inequalities 
    \begin{align}
    X_{ij} + X_{j \ell} \leq X_{i \ell} + x_j &\quad \forall \{i, j, \ell\} \subset C \label{triangle_inequalities_3} \\
    x_{i}+x_{j}+x_{\ell} \leq 1 + X_{i \ell} + X_{j \ell} + X_{ij} &\quad \forall \{i, j, \ell\} \subset C. \label{triangle_inequalities_4}
    \end{align}
    Then, the inequality~\eqref{stable_set_1} holds.
\end{lem}

\begin{proof}
Without loss of generality, let~$V = \{1, \ldots, n\}$, and let~$C = \{1, \dots, 5\}$, and~$E(C) = \{\{1,2\}, \{2,3\}, \{3,4\}, \{4,5\}, \{1,5\}\}$. Given that~$(x, X)$ is a feasible solution for the formulation~\eqref{theta_stable_set} strengthened with additional inequalities,~$X$ is symmetric,~$X_{ij} = 0$ for all edges~$\{i, j\} \in E(C)$, and~$x_i = X_{ii}$ for all~$i \in V$. Furthermore, since inequalities~\eqref{triangle_inequalities_3} and~\eqref{triangle_inequalities_4} are satisfied, we have that
\begin{align*}
X_{13} + X_{14} &\leq X_{11} \\
X_{24} + X_{25} &\leq X_{22} \\
X_{13} + X_{35} &\leq X_{33} \\
X_{11} + X_{22} + X_{33} &\leq 1 + X_{13}.
\end{align*}
Therefore,
\begin{align*}
\sum_{\substack{ i, j \in C \\ i < j}} X_{ij} + X_{13} \leq X_{11} + X_{22} + X_{33} \leq 1 + X_{13},
\end{align*}
so~\eqref{stable_set_1} holds. 
\end{proof}

The statement of Lemma~\ref{triangles_imply_sum_less_equal_one} demonstrates the power of triangle inequalities for obtaining strong upper bounds on the stability number of a graph. Nevertheless, following the results of Theorem~\ref{stability_number_2}, we can obtain even stronger bounds for~$5$-cycles. To show this, we first note the following direct consequence of Theorem~\ref{stability_number_2}.

\begin{cor}
Let~$G = (V, E)$ be a graph and let~$C \subseteq V$ induce a~$5$-cycle in~$G$. If~$X_C \in \STAB^2(G_C)$, then
	\begin{align}
		\sum_{i \in C} X_{ii} &\leq 1 + \sum_{ \{i, j\} \in \overline{E}(C)} X_{ij}. \label{cycles_new}
	\end{align}
\end{cor}

Now, it turns out that this inequality is not implied when we strengthen the Lovász theta function with triangle inequalities, as illustrated in the next example.

\begin{example} 
Let~$G$ be as in Example~\ref{first_example}. We set
    \begin{alignat*}{3}
		X_{ii} &= 0.3 &&\quad \forall i \in C \\
        x_i &= X_{ii} &&\quad \forall i \in C \\
		X_{ij} &= 0 &&\quad \forall \{i, j\} \in E(C) \\
		X_{ij} &= 0.075 &&\quad \forall \{i, j\} \in \overline{E}(C).
    \end{alignat*}
Then,~$\begin{pmatrix} X & x \\ x^\top & 1 \end{pmatrix} \succeq 0$,~$X_{ij} = 0$ for all~$\{i, j\} \in E(C)$, and~$\diag(X) = x$. Hence,~$(x, X)$ is feasible for~\eqref{theta_stable_set}. Furthermore, by construction, inequalities~\eqref{triangle_inequalities_1} and~\eqref{triangle_inequalities_2} imposed on all~$\{i, j, \ell\} \subset C$ are satisfied. However, 
\begin{align*}
    \sum_{i \in C} X_{ii} = 1.5 > 1.375 = 1 + \sum_{ \{i, j\} \in \overline{E}(C)} X_{ij},
\end{align*}
so~\eqref{cycles_new} is not satisfied.
\end{example}

As the next step, we consider subgraphs induced by an odd cycle~$C$ of arbitrary length and a vertex not contained in~$C$ and present two inequalities valid for such subgraphs. Moreover, we show that these inequalities are not implied by the Lovász theta function strengthened with triangle inequalities.

\begin{lem}\label{thm_stable_set_2}
	Let~$G = (V, E)$ be a graph, let~$C \subseteq V$ induce an odd cycle in~$G$, and let~$k \in V \setminus C$. Then, the inequalities 
	\begin{align}
		\sum_{i \in C} X_{ik} &\leq \frac{1}{2}(\vert C \vert - 1)X_{kk} \label{stable_set_2} \\
        \sum_{i \in C} X_{ii} + \frac{1}{2}(\vert C \vert - 1)X_{kk} &\leq \frac{1}{2}(\vert C \vert - 1) + \sum_{i \in C} X_{ik} \label{stable_set_3}
	\end{align}
    are valid for~$X_{C \cup \{k\}} \in \STAB^{2}(G_{C \cup \{k\}})$.
\end{lem}

\begin{proof}
    Without loss of generality, let~$V = \{ 1, \ldots, n \}$ and~$C = \{1, \ldots, k-1\}$. Let~$s = (s_1, \ldots, s_k)^\top$ be the incidence vector of a stable set in~$G_{C \cup \{k\}}$ and define~$X_{C \cup \{k\}} = ss^\top$. By construction, we have
    \begin{align*}
		\sum_{i \in C} X_{ik} = \sum_{i \in C} s_i s_k.
	\end{align*}

    Since~$C$ induces an odd cycle in~$G$, any stable set in~$G_C$ has at most~$\frac{1}{2}(\vert C \vert - 1)$ vertices. Hence,
    \begin{align}\label{odd_cycles_3}
		\sum_{i \in C} s_i \leq \frac{1}{2}(\vert C \vert - 1).
	\end{align}
    
    Multiplying this inequality by~$s_k$, we obtain 
    \begin{align}\label{odd_cycles_proof_1}
        \sum_{i \in C} s_i s_k \leq \frac{1}{2}(\vert C \vert - 1) s_k.
	\end{align}

    Now,~$s_j \in \{0, 1\}$ for any~$j \in C \cup \{k\}$, and therefore
     \begin{alignat*}{2}
    X_{ii} & = s_i s_i = s_i &\quad \forall i \in C \\
    X_{kk} &= s_k s_k = s_k.
    \end{alignat*}
    
    Hence, from~\eqref{odd_cycles_proof_1} we obtain
    \begin{align*}
    \sum_{i \in C} X_{ik} &\leq \frac{1}{2}(\vert C \vert - 1)X_{kk}.
    \end{align*}    
    Thus,~\eqref{stable_set_2} holds. 
    
    To obtain~\eqref{stable_set_3}, we proceed in the same manner and multiply~\eqref{odd_cycles_3} by~$(1 - s_k)$, yielding
    \begin{align*}
		\sum_{i \in C} s_i + \frac{1}{2}(\vert C \vert - 1) s_k \leq \frac{1}{2}(\vert C \vert - 1) + \sum_{i \in C} s_i s_k.
	\end{align*}
    
    Hence,
    \begin{align*}
    \sum_{i \in C} X_{ii} + \frac{1}{2}(\vert C \vert - 1)X_{kk} &\leq \frac{1}{2}(\vert C \vert - 1) + \sum_{i \in C} X_{ik},
    \end{align*}
    so~\eqref{stable_set_3} also holds.

    Finally, since~\eqref{stable_set_2} and~\eqref{stable_set_3} hold for all rank-one matrices~$X_{C \cup \{k\}} = ss^\top$ corresponding to stable set matrices in~$\mathcal{S}^2(G_{C \cup \{k\}})$, they remain valid also for all~$X_{C \cup \{k\}} \in \STAB^{2}(G_{C \cup \{k\}})$.
\end{proof}

To show that the inequalities~\eqref{stable_set_2} and~\eqref{stable_set_3} are not implied by the Lovász theta function strengthened with triangle inequalities, we give two counterexamples.

\begin{example}\label{counterexample_one}
Let~$G$ be a graph with~$V = \{1, \ldots, 6\}$, let~$C \subseteq V$ be a~$5$-cycle with~$C = \{1, \dots, 5\}$, and let~$E(G) = E(C) = \{\{1,2\}, \{2,3\}, \{3,4\}, \{4,5\}, \{1,5\}\}$. We define
\begin{alignat*}{3}
	X_{ii} &= 0.3333 &&\quad \forall i \in C \\
	X_{ij} &= 0 &&\quad \forall \{i, j\} \in E(C) \\
	X_{ij} &= 0.1667 &&\quad \forall \{i, j\} \in \overline{E}(C) \\
        X_{66} &= 0.4444 && \\
	X_{i6} = X_{6i} &= 0.2222 &&\quad \forall i \in C \\
    x_i &= X_{ii} &&\quad \forall i \in V.
\end{alignat*}
Then,~$\begin{pmatrix} X & x \\ x^\top & 1 \end{pmatrix} \succeq 0$,~$X_{ij} = 0$ for all~$\{i, j\} \in E(G)$, and~$\diag(X) = x$. Furthermore, inequalities~\eqref{triangle_inequalities_1} and~\eqref{triangle_inequalities_2} are satisfied. Hence,~$(x, X)$ is feasible for the formulation~\eqref{theta_stable_set} strengthened with triangle inequalities. However,
\begin{align*}
\sum_{i \in C} X_{i 6} = 1.1110 > 0.8888 = 2X_{66}.
\end{align*}
Thus,~\eqref{stable_set_2} is not satisfied.
\end{example}

\begin{example}\label{counterexample_two}
Let~$G$ be as in Example~\ref{counterexample_one} and let
\begin{alignat*}{3}
	X_{ii} &= 0.3333 &&\quad \forall i \in C \\
	X_{ij} &= 0 &&\quad \forall \{i, j\} \in E(C) \\
	X_{ij} & = 0.1667 &&\quad \forall \{i, j\} \in \overline{E}(C) \\
        X_{66} &= 0.5556 && \\
	X_{i6} = X_{6i}&= 0.1111 &&\quad \forall i \in C \\
        x_i &= X_{ii} &&\quad \forall i \in V.
\end{alignat*}
Then, it can be easily shown that~$(x, X)$ is feasible for the formulation~\eqref{theta_stable_set} strengthened with inequalities~\eqref{triangle_inequalities_1} and~\eqref{triangle_inequalities_2}. However,
\begin{align*}
\sum_{i \in C} X_{i i} + 2X_{66} = 2.7777 > 2.5555 = 2 + \sum_{i \in C} X_{i 6}.
\end{align*}
Therefore,~\eqref{stable_set_2} does not hold.
\end{example}

The next question that naturally arises is what happens when we consider an odd cycle of arbitrary length and a vertex~$k$ that is an element of~$C$. The answer is given in the next statement. 

\begin{lem}\label{lemma_odd_cycles}
	Let~$G = (V, E)$ be a graph, let~$C \subseteq V$ induce an odd cycle in~$G$, and let~$k \in C$. Then, the inequalities 
	\begin{align}
		\sum_{\substack{i \in C \\ i \notin N[k]}} X_{ik} \leq \frac{1}{2}(\vert C \vert - 3) X_{kk} \label{stable_set_4} \\
        \sum_{i \in C} X_{ii} + \frac{1}{2}(\vert C \vert - 3)X_{kk} &\leq \frac{1}{2}(\vert C \vert - 1) + \sum_{\substack{i \in C \\ i \notin N[k]}} X_{ik} \label{stable_set_5}
	\end{align}
    are valid for~$X_C \in \STAB^{2}(G_{C})$.
\end{lem}

\begin{proof}
    Without loss of generality, let~$V = \{ 1, \ldots, n \}$ and~$C = \{1, \ldots, k\}$. Let~$s = (s_1, \ldots, s_k)^\top$ be the incidence vector of a stable set in~$G_{C}$ and define~$X_{C} = ss^\top$. Then, by construction,
    \begin{align}\label{vertex_in_cycle}
        \sum_{i \in C} X_{ik} = \sum_{i \notin N[k]} X_{ik} + X_{kk} = \sum_{i \in C} s_i s_k.
    \end{align}

    Now, given that~$C$ is an odd cycle,~\eqref{odd_cycles_3} also holds. Using a similar argument as in the proof of Lemma~\ref{thm_stable_set_2}, we multiply~\eqref{odd_cycles_3} by~$s_k$ and obtain 
    \begin{align*}
        \sum_{i \in C} s_i s_k \leq \frac{1}{2}(\vert C \vert - 1) s_k.
    \end{align*}

    Now, as argued in the proof of Lemma~\ref{thm_stable_set_2}, we have
    \begin{alignat*}{2}
    X_{ii} & = s_i s_i = s_i &\quad \forall i \in C \\
    X_{kk} &= s_k s_k = s_k.
    \end{alignat*}
    
    Hence, incorporating~\eqref{vertex_in_cycle} and rearranging the terms, we arrive at inequality~\eqref{stable_set_4}. 
    
    Next, we multiply~\eqref{odd_cycles_3} by~$(1-s_k)$, and after rearranging and using the previous arguments, we obtain inequality~\eqref{stable_set_5}. 

    Finally, given that~\eqref{stable_set_4} and~\eqref{stable_set_5} hold for all rank-one matrices~$X_C = ss^\top$ corresponding to stable set matrices in~$\mathcal{S}^2(G_C)$, they are valid for all~$X_C \in \STAB^2(G_C)$, which completes the proof.
\end{proof}

Hence, Lemma~\ref{lemma_odd_cycles} gives an insight into what happens when we consider an odd cycle of arbitrary length and a vertex that is an element of~$C$. However, it turns out that adding triangle inequalities into the SDP formulation~\eqref{theta_stable_set} implies the inequalities presented in Lemma~\ref{lemma_odd_cycles}.

\begin{lem}\label{lemma_odd_cycles_2}
	Let~$G = (V, E)$ be a graph, let~$C \subseteq V$ induce an odd cycle in~$G$, let~$j \in C$, and let~$(x, X)$ be a feasible solution for~\eqref{theta_stable_set} strengthened with inequalities~\eqref{triangle_inequalities_3} and~\eqref{triangle_inequalities_4}. Then, the inequalities~\eqref{stable_set_4} and~\eqref{stable_set_5} are satisfied.
\end{lem}

\begin{proof}
Without loss of generality, let~$V = \{1, \ldots, n\}$,~$C = \{1, \ldots, k\}$, and let 
\begin{align*}E(C) = \left\{ \{s, t\} \subset C : t = s + 1 ~\forall s \in C \setminus \{k\} \right\} \cup \{ \{1, k\} \}.
\end{align*}
Furthermore, let~$j \in C$ be an arbitrary but fixed vertex in~$C$. 

We define a subset of edges~$E' \subset E(C)$ such that no edge in~$E'$ contains any vertex from~$N[j]$, such that none of the edges in~$E'$ share a common vertex, and such that~$E'$ is as large as possible, i.e.,~$\vert E' \vert = \frac{1}{2}(\vert C \vert - 3)$. For instance, if~$j = 1$, then~$E' = \left\{  \{3, 4\}, \{5, 6\}, \ldots, \{k-2, k-1\} \right\}$. 

Now, since $\diag(X) = x$, summing up the imposed triangle inequalities~\eqref{triangle_inequalities_3} over the set~$E'$ yields
\begin{align}\label{inequality_sum_cycles}
\sum_{\{i, \ell\} \in E'} (X_{ij} + X_{j \ell}) \leq \sum_{\{i ,\ell\} \in E'} (X_{i\ell} + X_{jj}).
\end{align}

Given that~$(x, X)$ is a feasible solution for~\eqref{theta_stable_set}, we have that~$X_{i \ell} = 0$ for all edges~$\{i, \ell\} \in E'$. Furthermore,~$\vert E' \vert = \frac{1}{2}(\vert C \vert - 3)$ according to construction. Hence, we can rewrite~\eqref{inequality_sum_cycles} as
\begin{align*}
\sum_{\{i, \ell\} \in E'} (X_{ij} + X_{j \ell}) \leq \frac{1}{2}(\vert C \vert - 3) X_{jj}.
\end{align*}

Now, regarding the left-hand side of this inequality, we note that each edge~$\{i, \ell\} \in E'$ is chosen such that~$i, \ell \notin \text{N}[j]$ and no two edges in~$E'$ share a common vertex. Given this construction, each term~$X_{ij} + X_{j \ell}$ in the sum over~$E'$ captures the contribution from pairs~$\{i, j\}$ and~$\{j, \ell\}$, where~$i$ and~$\ell$ are distinct vertices in~$C \setminus N[j]$.

Since each edge in~$E'$ connects two distinct vertices in~$C \setminus N[j]$, summing up over all edges in~$E'$ ensures that each vertex~$i \in C \setminus N[j]$ appears exactly once in the form~$X_{ij}$. Therefore, 
\begin{align*}
\sum_{\{i, \ell\} \in E'} (X_{ij} + X_{j \ell}) = \sum_{\substack{i \in C \\ i \notin N[j]}} X_{ij}.  
\end{align*}

Hence, altogether, we obtain 
\begin{align*}
\sum_{\substack{i \in C \\ i \notin N[j]}} X_{ij} \leq \frac{1}{2}(\vert C \vert - 3) X_{jj}. 
\end{align*}
Thus,~\eqref{stable_set_4} holds.

To prove that~\eqref{stable_set_5} also holds, we again consider the vertex~$j$ and define a subset of edges~$\tilde{E} \subset E(C)$ such that no edge in~$\tilde{E}$ contains the vertex~$j$, such that none of the edges in~$\tilde{E}$ share a common vertex, and such that~$\tilde{E}$ is as large as possible, i.e.,~$\vert \tilde{E} \vert = \frac{1}{2}(\vert C \vert - 1)$. For instance, if~$j = 1$, then~$\tilde{E} = \left\{  \{2, 3\}, \{4, 5\}, \ldots, \{k-1, k\} \right\}$.

Analogously to the previous argumentation, we sum up imposed triangle inequalities~\eqref{triangle_inequalities_4} over the set~$\tilde{E}$ and obtain
\begin{align}\label{triangle_cycles_2}
\sum_{\{i, \ell\} \in \tilde{E}} (X_{ii}+X_{jj}+X_{\ell \ell}) \leq \sum_{\{i, \ell\} \in \tilde{E}} (1 + X_{i \ell} + X_{j \ell} + X_{i j}). 
\end{align}

Now, we have that~$X_{i \ell} = 0$ for all edges in~$\tilde{E}$. Moreover,~$X_{j \ell} = 0$ or~$X_{i j} = 0$ if~$i$ or~$\ell$ are neighbors of the vertex~$j$, respectively. Furthermore, according to the construction, there are altogether~$\frac{1}{2}(\vert C \vert - 1)$ edges in the set~$\tilde{E}$. 

Regarding the term~$X_{j \ell} + X_{ij}$ in the sum over~$\tilde{E}$, from the construction, we have that each vertex~$i \in C \setminus \{j\}$ appears exactly once in the form~$X_{ij}$. Thus, we can express~\eqref{triangle_cycles_2} as
\begin{align*}
\sum_{\{i, \ell\} \in \tilde{E}} (X_{ii}+X_{jj}+X_{\ell \ell}) \leq \frac{1}{2}(\vert C \vert- 1) + \sum_{\substack{i \in C \\ i \notin N[j]}} X_{ij}.
\end{align*}

Finally, following the same argument as above, we have that
\begin{align*}
\sum_{\{i, \ell\} \in \tilde{E}} (X_{ii}+X_{jj}+X_{\ell \ell}) = \sum_{\substack{i \in C \\ i \neq j}} X_{ii} + \frac{1}{2}(\vert C \vert - 1) X_{jj}.
\end{align*}

Therefore, we obtain
\begin{align*}
\sum_{i \in C} X_{ii} + \frac{1}{2}(\vert C \vert - 3) X_{jj} \leq \frac{1}{2}(\vert C \vert - 1) + \sum_{\substack{i \in C \\ i \notin N[j]}} X_{ij},
\end{align*}
which finishes the proof.
\end{proof}

\subsection{Subgraphs Containing Odd Antiholes}\label{section_antiholes}

Subgraphs induced by odd antiholes are another class of subgraphs with stability number two. Thus, we can use Theorem~\ref{stability_number_2} to completely describe their convex hulls. Furthermore, by Lemma~\ref{thm_stable_set_1}, we know that the inequality~\eqref{stable_set_inequality} is valid for such subgraphs.

However, if we consider a subgraph induced by an odd antihole along with an additional vertex~$k$ not contained in the antihole, the stability number of this induced subgraph is no longer guaranteed to be two. We apply the same approach as for odd cycles in the previous section and obtain the following result.

\begin{lem}\label{lemma_antihole}
	Let~$G = (V, E)$ be a graph, let~$A \subseteq V$ induce an odd antihole in~$G$, and let~$k \in V \setminus A$. Then, the following inequalities 
	\begin{align}
		\sum_{i \in A} X_{ik} &\leq 2 X_{kk} \label{antiholes_1} \\
        \sum_{i \in A} X_{ii} + 2 X_{kk} &\leq 2 + \sum_{i \in A} X_{ik} \label{antiholes_2}
	\end{align}
    are valid for~$X_{A \cup \{k\}} \in \STAB^{2}(G_{A \cup \{k\}})$.
\end{lem}

\begin{proof}
   Without loss of generality, let~$V = \{ 1, \ldots, n \}$ and~$A = \{1, \ldots, k-1\}$. Let~$s = (s_1, \ldots, s_k)^\top$ be the incidence vector of a stable set in~$G_{A \cup \{k\}}$ and define~$X_{A \cup \{k\}} = ss^\top$. Then, by construction,
    \begin{align*}
		\sum_{i \in A} X_{ik} = \sum_{i \in A} s_i s_k.
	\end{align*}

    Given that~$A$ induces an odd antihole in~$G$, any stable set in~$G_A$ has at most~$2$ vertices. Thus,
    \begin{align}\label{odd_antiholes_3}
		\sum_{i \in C} s_i \leq 2.
	\end{align}
    
    Following the same steps and argumentation as in the proof of Lemma~\ref{thm_stable_set_2}, we multiply~\eqref{odd_antiholes_3} by~$s_k$ and $(1 - s_k)$ and obtain that~\eqref{antiholes_1} and~\eqref{antiholes_2} are valid for the matrix $X_{A \cup \{k\}}$.  

    Hence,~\eqref{antiholes_1} and~\eqref{antiholes_2} are valid for all rank-one matrices~$X_{A \cup \{k\}} = ss^\top$ corresponding to stable set matrices in~$\mathcal{S}^2(G_{A \cup \{k\}})$. Thus, they are valid for all~$X_{A \cup \{k\}} \in \STAB^2(G_{A \cup \{k\}})$, completing the proof. 
\end{proof}

Note that a complement of a~$5$-hole is again a~$5$-hole. Therefore, from Examples~\ref{counterexample_one} and~\ref{counterexample_two} we know that inequalities~\eqref{antiholes_1} and~\eqref{antiholes_2} are not implied by triangle inequalities. Hence, adding the inequalities~\eqref{antiholes_1} and~\eqref{antiholes_2} into the SDP to compute the Lovász theta function might yield a tighter upper bound on~$\alpha(G)$.

Finally, recall that for odd cycles, we also examined the case where the vertex~$k$ is part of the cycle. For odd antiholes, however, this is unnecessary, because the underlying subgraph has stability number two, for which the inequalities from Theorem~\ref{stability_number_2} already fully describe the convex hull.

\section{Valid Inequalities for the Graph Coloring Problem}\label{inequalities_coloring}

We now turn our attention to the graph coloring problem and investigate subgraphs containing odd cycles, proposing two inequalities valid for this problem. Similar to the approach taken for the stable set problem, the first inequality we discuss is valid for odd cycles of length~$5$, while the second holds for odd cycles of arbitrary length along with an additional vertex that is not contained in the cycle.

\begin{lem}\label{thm_coloring_1}
	Let~$G = (V, E)$ be a graph and let~$C \subseteq V$ induce a~$5$-cycle in~$G$. Then the inequality 
	\begin{align}
		\sum_{\substack{i,j \in C \\ i < j}} X_{ij} \leq 2 \label{coloring_1}
	\end{align}
    is valid for~$X_C \in \COL(G_C)$.
\end{lem}

\begin{proof}  
Let~$X$ be a coloring matrix. By definition, we can express~$X$ as  
\begin{align}\label{matrix_X_coloring}  
X = \sum_{p = 1}^{\ell} s_p s_p^\top,  
\end{align}  
where each~$s_p$ is the incidence vector of a stable set in~$G$, with~$s_p \neq 0$ and~$\sum_{p = 1}^{\ell} s_p = \mathds{1}_n$.  

Hence, the sum  
\begin{align*}  
\sum_{\substack{i,j \in C \\ i < j}} X_{ij}  
\end{align*}  
counts the number of indices~$p$ for which the set~$s_p$ contains at least two vertices from~$C$.

Since the stability number of~$C$ is at most 2, and~$C$ is a~$5$-cycle, there can be at most two indices~$p$ for which this condition holds. Therefore, we conclude that~\eqref{coloring_1} is valid for any coloring matrix~$X_C \in \mathcal{C}(G_C)$, and therewith for any~$X_C \in \COL(G_C)$, which completes the proof.
\end{proof}  

\begin{lem}\label{thm_coloring_2}
Let~$G = (V, E)$ be a graph, let~$C \subseteq V$ induce an odd cycle in~$G$, and let~$k \in V$ such that~$k \notin C$. Then the inequality
	\begin{align}
		\sum_{i \in C} X_{ik} \leq \frac{1}{2} (\vert C \vert - 1).\label{coloring_2}
	\end{align}
    is valid for~$\COL(G_{C \cup \{k\}})$.
\end{lem}

\begin{proof}  
Let~$X$ be a coloring matrix. By definition, and as noted in the proof of Theorem~\ref{thm_coloring_1}, we can express~$X$ as~\eqref{matrix_X_coloring}, where each~$s_p$ is the incidence vector of a stable set in~$G$, with~$s_p \neq 0$ and~$\sum_{p = 1}^{\ell} s_p = \mathds{1}_n$.  

Since~$k$ is a vertex in the graph, it must belong to at least one stable set. That is, there exists some index~$p \in \{1, \ldots, \ell\}$ such that~$k \in s_p$. Without loss of generality, assume that~$k \in s_1$. Then, we can express the sum  
\begin{align*}  
\sum_{i \in C} X_{ik}  
\end{align*}  
as the number of vertices in~$C$ that belong to the stable set~$s_1$.

Given that~$C$ is an odd cycle, its stability number is at most~$\frac{1}{2} (\vert C \vert - 1)$. Thus, the number of vertices in~$C$ that can belong to the same stable set as~$k$ is at most~$\frac{1}{2} (\vert C \vert - 1)$. This directly implies that  
\begin{align*}  
\sum_{i \in C} X_{ik} \leq \frac{1}{2} (\vert C \vert - 1).
\end{align*}  

Therefore,~\eqref{coloring_2} holds for any coloring matrix~$X_C \in \mathcal{C}(G_C)$. Thus,~\eqref{coloring_2} is valid for any~$X_C \in \COL(G_{C \cup \{k\}})$.
\end{proof}  

In the next section, we will incorporate the proposed inequalities \eqref{clique_vertex_four}, \eqref{coloring_1}, and \eqref{coloring_2} into the SDP to compute the Lovász theta function and evaluate the resulting bounds. Before doing so, we will first strengthen the formulation~\eqref{theta_coloring} by adding the previously known non-negativity constraints \eqref{non-n-col} as well as the triangle inequalities \eqref{triangle_inequalities_coloring}, as these can be incorporated efficiently.

In the previous section, we showed that some of the proposed inequalities for the stable set problem hold when the SDP for computing the Lovász theta function~\eqref{theta_stable_set} is strengthened with triangle inequalities. This motivated a careful examination of some of the proposed inequalities to determine whether they are implied by the triangle inequalities. However, this relationship does not hold for the graph coloring problem, i.e., for the inequalities introduced in this section. For brevity, we do not provide explicit examples, but our computational study will confirm this observation.    
    
\section{Computational Study}\label{Computational_results}

We now perform a computational study to evaluate the strengthening of the Lovász theta function towards the stability number $\alpha(G)$ and chromatic number $\chi(G)$ by adding valid inequalities into the SDP formulation used for its computation. The computations are carried out on an Intel(R) Core(TM) i7-1260P CPU @ 2.10GHz with 64GB of RAM. We use Matlab for implementation, with MOSEK serving as the SDP solver. 

\subsection{The Stable Set Problem}\label{results_stable_set}

Our computational methodology is structured as follows. We strengthen the Lovász theta function~$\vartheta(G)$ towards the stability number~$\alpha(G)$ by adding several classes of inequalities to the SDP for its computation~\eqref{theta_stable_set}. Specifically, we add the non-negativity constraints~\eqref{non-n}, the triangle inequalities~\eqref{triangle_inequalities_1} and \eqref{triangle_inequalities_2}, and the proposed inequalities~\eqref{cliques_3}, \eqref{cliques_final}, \eqref{cycles_new}, \eqref{stable_set_2}, and \eqref{stable_set_3}. To incorporate these inequalities efficiently, we enumerate all cliques with at most five vertices and consider only chordless cycles of length~$5$. Since an odd antihole of length~$5$ is a chordless~$5$-cycle, we do not consider any inequalities proposed in Section~\ref{section_antiholes}. The addition of inequalities is done in two phases:

\begin{itemize}
\item 
%\noindent
\textbf{Phase 1:} In this phase, we incorporate well-known inequalities that have been investigated by other researchers. After computing the optimal solution for the formulation~\eqref{theta_stable_set}, we examine the obtained solution and search for violations of the non-negativity constraints~\eqref{non-n} as well as the triangle inequalities~\eqref{triangle_inequalities_1} and~\eqref{triangle_inequalities_2}. Since the number of found violations can be quite large, we add violated inequalities iteratively to the SDP~\eqref{theta_stable_set}. Specifically, we add all violations of~\eqref{non-n} in each iteration. For~\eqref{triangle_inequalities_1} and~\eqref{triangle_inequalities_2}, we consider violations greater than~$0.025$, and we add up to~$2n$ of the largest violations per inequality type and per iteration. We iterate until there are fewer than~$n$ violations of all considered inequalities or until the number of iterations reaches~$10$. The obtained upper bound on~$\alpha(G)$ should satisfy almost all inequalities~\eqref{non-n}, \eqref{triangle_inequalities_1}, and~\eqref{triangle_inequalities_2}. We denote this bound by BOUND 1.
\item
%\noindent
\textbf{Phase 2:} We now add the proposed inequalities~\eqref{cliques_3} and~\eqref{cliques_final}, \eqref{cycles_new}, \eqref{stable_set_2} and~\eqref{stable_set_3}. As already mentioned, to separate these inequalities efficiently, we enumerate all cliques with at most five vertices and consider only chordless cycles of length~$5$. 
Hence, we start from the final solution of the SDP from Phase 1 and examine violations of the mentioned inequalities. Then, similar to the procedure from Phase 1, we consider violations of these inequalities greater than~$0.025$, and we add up to~$2n$ of the largest violations per inequality type and iteration. This procedure is repeated until there are fewer than~$n$ violations of all considered inequalities or until the number of iterations reaches~$10$. We denote the resulting upper bound on~$\alpha(G)$ by BOUND 2.
\end{itemize}

\subsubsection{First Set of Experiments}

In our first set of experiments, we deal with some instances considered in~\cite{Gaa:20}:
\begin{itemize}
\item Near-regular graphs: These graphs were generated by the following procedure. First, a perfect matching on~$nr$ vertices was selected. Then, consecutive groups of~$r$ vertices were combined into a single vertex, which resulted in a regular multigraph on~$n$ vertices. Finally, to obtain near-regular graphs, loops and multiple edges from this multigraph were removed.
\item Random graphs from the Erdős–Rényi model~$G(n, p)$: In this model, the number of vertices~$n$ is fixed, and each edge is included with probability~$p$ independently of all other edges. 
\item Torus graphs: These instances are constructed as follows. For a given~$d$, the graph~$T_d$ has~$d^2$ vertices labeled by~$(i, j)$ for~$i, j \in \{1, \ldots, d\}$. The vertical edges connect vertices with neighboring~$i$ indices (keeping~$j$ fixed), creating edges~$\{(i, j), (i+1, j)\}$ modulo~$d$. Similarly, the horizontal edges connect vertices with neighboring~$j$ indices (keeping~$i$ fixed), resulting in edges of the form~$\{(i, j), (i, j + 1)\}$ modulo~$d$. Thus, the graph has~$n = d^2$ vertices and~$m = 2n$ edges. As noted in~\cite{Gaa:20}, for even~$d$,~$\vartheta(T_d) = \alpha(T_d)$. Therefore, the study in~\cite{Gaa:20} focused on torus graphs with odd~$d$. 
\end{itemize}

The goal of this computational study is to investigate the quality of obtained bounds BOUND 1 and BOUND 2 as upper bounds on~$\alpha(G)$ and to compare them to those presented in~\cite{Gaa:20}, which we denote by GR. Specifically, we focus on two key aspects: first, the improvement of BOUND 2 over BOUND 1, given that BOUND 2 is a refinement of BOUND 1; and second, the comparison between BOUND 2 and GR. Additionally, we assess the computational times required for these bounds. Although our computational setup differs from that in~\cite{Gaa:20}, we include the reported computational times from~\cite{Gaa:20} to provide context for the computational efforts involved.

The results for this first set of experiments are presented in Table~\ref{table_1}. The first columns of Table~\ref{table_1} provide general information about the instances, including the name of the instance, the number of vertices~$n$, the number of edges~$m$, the value of~$\alpha(G)$ as given in \cite{Gaa:20}, and the value of the Lovász theta function~$\vartheta(G)$. Next, we report the computed bounds BOUND 1 and BOUND 2, along with the bounds GR from~\cite{Gaa:20}. For each instance, the best computed bound on~$\alpha(G)$ is highlighted in bold. Additionally, if BOUND 2 yields a better integer bound on~$\alpha(G)$ than bound GR, the corresponding cell is shaded in gray. Finally, we present the computational times (time) required for the computation of bounds BOUND 1 and BOUND 2, as well as the times for the computation of the GR bound, as reported in \cite{Gaa:20}.

    \begin{table}[h!]
	\caption{Stable set results for instances considered in~\cite{Gaa:20}}
	\footnotesize
	\begin{center}
        %\begin{adjustbox}{width=\textwidth}
		\begin{tabular}{|l r r c c | rr| rr|rr|}
		\hline
    Graph &~$n$ &~$m$ &~$\alpha(G) \geq$ &~$\vartheta(G)$  & BOUND 1&(time)& BOUND 2&(time) & GR& (time)\\
		\hline
    reg\_n100\_r4 & 100 & 195 & 40 & 43.449  & 41.246 & (17) & 40.753 & (22) & \textbf{40.687} & (1164) \\
        \hline
    reg\_n100\_r6 & 100 & 294 & 34 & 37.815 & 36.224 & (7) & \textbf{35.046} & (26) & 35.246 & (1062)\\
        \hline
    reg\_n100\_r8 & 100 & 377 & 31 & 34.480 & 33.337& (4) & \textbf{32.031} & (22) & 32.190 & (1067) \\  
        \hline
    reg\_n200\_r4 & 200 & 400 & 80 & 87.759 & 83.498 & (111) & \cellcolor{lightgray}{\textbf{82.292}} & (196) & 83.772& (1437) \\
        \hline
    reg\_n200\_r6 & 200 & 593 & 68 & 79.276 & 76.047& (25) & \cellcolor{lightgray}{ \textbf{73.738}} & (132) & 75.555 & (1523) \\
        \hline
    reg\_n200\_r8 & 200 & 792 & 60 & 70.790 & 69.110 & (11) & \cellcolor{lightgray}{\textbf{66.743}} & (83) & 67.785 & (1944) \\
        \hline
    reg\_n200\_r10& 200 & 980 & 57 & 66.418 & 65.142 & (6) & \textbf{62.714} & (56) & 62.894 & (2556) \\
        \hline
        rand\_n100\_p004 & 100 & 212 & 45 & 46.067 & 45.032 & (22) & \textbf{45.006} & (6) & 45.021 & (432) \\
        \hline
        rand\_n100\_p006 & 100 & 303 & 38 & 40.361 & 38.909 & (13) & \textbf{38.435} & (14) & 38.439 & (887) \\
        \hline
        rand\_n100\_p008 & 100 & 443 & 32 & 34.847 & 33.575 & (5) & \textbf{32.431} & (18) & 32.579 & (1262) \\
        \hline
        rand\_n100\_p010 & 100 & 489 & 32 & 34.020 & 32.934 & (5) & \textbf{32.126} & (19) & 32.191 & (1138) \\
        \hline
        rand\_n200\_p002 & 200 & 407 & 95 & 95.778 & 95.044 & (222) & 95.044 & (1) & \textbf{95.032} & (836) \\
        \hline
        rand\_n200\_p003 & 200 & 631 & 80 & 83.662 & 81.560 & (39) &  \textbf{81.064} & (57) & 81.224 & (1867) \\
        \hline
        rand\_n200\_p004 & 200 & 816 & 67 & 73.908 & 71.654 & (17) &  \cellcolor{lightgray} \textbf{69.854} & (79) & 70.839 & (2227) \\
        \hline
        rand\_n200\_p005 & 200 & 991 & 62 & 69.039 & 67.313 & (19) &  \cellcolor{lightgray} \textbf{65.509} & (82) & 66.091 & (2411) \\
        \hline        
         torus\_5 & 25 & 50 & 10 & 11.180 & \textbf{10.000} & (1) &
        \textbf{10.000}  & (1) & 10.002 & (33) \\
        \hline
        torus\_7 & 49 & 98 & 21 & 23.224 & \textbf{21.000} & (2)& 
         \textbf{21.000}  & (1) & 21.009 & (127) \\
        \hline
        torus\_9 & 81 & 162 & 36 & 39.241 & \textbf{36.000} & (24) &
         \textbf{36.000} & (7) & 36.021 & (344) \\
        \hline
        torus\_11 & 121 & 242 & 55 & 59.249 & 55.022 & (81) &
        \textbf{55.019} & (19) & 55.066 & (851)\\
        \hline
        torus\_13 & 169 & 338 & 78 & 83.254 & 78.379 & (337) &
       \cellcolor{lightgray} \textbf{78.048} & (129) & 79.084 & (1031) \\
        \hline
        torus\_15 & 225 & 450 & 105 & 111.257 & 108.208 & (1517) &
        \cellcolor{lightgray} \textbf{105.214} & (504)  & 106.287 &(1615)\\
        \hline
		\end{tabular}
        %\end{adjustbox}
	\end{center}
	\label{table_1}
    \end{table}

    From the results for near-regular graphs presented in Table~\ref{table_1}, it is evident that BOUND 2 provides significantly tighter estimates on~$\alpha(G)$ compared to BOUND 1 across all tested instances. Notably, BOUND 2 outperforms BOUND 1 by at least one integer for every instance considered, while maintaining a generally low computational cost. Additionally, BOUND 2 frequently outperforms the GR bound, delivering better bounds in six out of seven cases. Even when BOUND 2 and the GR bound are equal, both bounds yield the same integer upper bound, which coincides with~$\alpha(G)$. Furthermore, BOUND 2 is more computationally efficient than the GR method, making it a strong choice for improving upper bounds on the stability number of near-regular graphs.

    The results for random graphs mirror those for near-regular graphs, i.e., BOUND 2 consistently provides tighter estimates of~$\alpha(G)$ compared to BOUND 1. Additionally, BOUND 2 frequently outperforms the GR bound, yielding better upper bounds on~$\alpha(G)$ in seven out of eight instances. Even in cases where GR performs better, both methods yield the same integer bound, matching~$\alpha(G)$. Additionally, BOUND 2 offers a significant computational advantage, effectively delivering tighter upper bounds on the stability numbers of random graphs while reducing computational effort.
   
    Finally, the results for two-dimensional torus graphs presented in Table~\ref{table_1} show that BOUND 2 consistently provides the tightest upper bounds on~$\alpha(G)$ across all six instances, outperforming the GR bound and delivering integer bounds that exactly match~$\alpha(G)$ for all evaluated torus graphs. Interestingly, BOUND 1 also proves to be highly competitive for torus graphs, often offering bounds that are better or very close to GR ones while requiring significantly less computational effort.

    In~\cite{Gaa:20}, the authors emphasized the high computational demands of calculating upper bounds on the stability numbers of torus graphs, suggesting the need for a specialized SDP solver. Our findings indicate that incorporating non-negativity constraints and triangle inequalities into the SDP for the Lovász theta function already improves the precision of these upper bounds while reducing computational time. Therefore, any specialized SDP solver developed for these instances should integrate these inequalities to enhance both accuracy and efficiency.

    \subsubsection{Second Set of Experiments}
    
    In our second set of experiments, we investigate several well-established instances from the literature:
    \begin{itemize}
    \item DIMACS Instances: We analyze several benchmarks from the Second DIMACS Implementation Challenge~\cite{Johnson1996}, held in 1992 and 1993, which focused on NP-hard problems including maximum clique, graph coloring, and satisfiability. For our study on the stable set problem, we use the complements of the maximum clique instances from this challenge.
    \item Spin graphs: The stability numbers of certain spin graphs were recently analyzed using an SDP-based branch-and-bound framework in~\cite{Gaar2022}. These instances are particularly interesting because, for most of the spin graphs studied in~\cite{Gaar2022}, only lower and upper bounds on~$\alpha(G)$ were provided. Therefore, it is valuable to investigate whether the upper bounds on~$\alpha(G)$ can be improved by strengthening the Lovász theta function~$\vartheta(G)$ through the addition of the proposed inequalities.
    \item Evil graphs: We consider extremely hard and versatile (evil) instances for the maximum clique problem recently introduced by Szabo and Zaválnij in~\cite{Szabo}, known for their wide range of difficulties. These instances, comprising a total of~$40$ graphs, are highly valued in algorithmic research. The clique numbers of evil graphs were investigated by San Segundo, Furini, {\'A}lvarez, and Pardalos in~\cite{SanSegundo2022}. Since our focus is on the stable set problem, we consider the complements of some of these instances. To the best of our knowledge, this is the first computational study of the SDP-based upper bounds on the stability numbers of evil instances.
    \end{itemize}

    The primary objective of this computational study is to evaluate the extent to which BOUND 2 improves upon BOUND 1 as an upper bound on~$\alpha(G)$. Additionally, we examine the computational times required to calculate both BOUND 1 and BOUND 2. The computational results for the spin graphs and selected DIMACS instances are presented in Table~\ref{table_4}. The table follows the same structure as those from the first set of experiments, with the exception that the GR bound is not reported. Additionally, we know that BOUND 2 outperforms BOUND 1 as an upper bound on~$\alpha(G)$ and that is why we do not use bold formatting to highlight the best computed bound for each instance. However, if BOUND 2 surpasses BOUND 1 by at least one integer value, we highlight this result in gray. 

    \begin{table}[h!]
	\caption{Stable set results for selected instances from the literature}
	\footnotesize
	\begin{center}
        %\begin{adjustbox}{width=\textwidth}
	\begin{tabular}{|l r r c r | rr| rr|}
	\hline
    Graph &~$n$ &~$m$ &~$\alpha(G)$ &~$\vartheta(G)$ &
    BOUND 1&(time) & BOUND 2& (time)\\
    \hline
    spin5 & 125 & 375 & 50 & 55.902 & 50.000 & (17) & 50.000&(6) \\
    \hline
    spin7 & 343 & 1029 &~$147 \leq \alpha \leq 151$ & 162.566& 147.000 & (1225) & 147.000 & (639) \\
    \hline
    MANN\_a9 & 45 & 72 & 16 & 17.475 & 17.220 & (2) & 17.220 & (1) \\
    \hline
    MANN\_a27 & 378 & 702 & 126 & 132.763 & 131.709 & (781) & 131.112 & (874) \\
    \hline
    C125.9 & 125 & 787 & 34 & 37.805 & 36.920 & (4) & \cellcolor{lightgray} 35.542 & (41) \\ 
    \hline
    C250.9 & 250 & 3141 & 44 & 56.241 & 55.771 & (18) & \cellcolor{lightgray} 54.899 & (479) \\
    \hline
    sanr200\_0\_9 & 200 & 2037 & 42 & 49.274& 48.723 & (11) &\cellcolor{lightgray} 47.474 & (229) \\
    \hline
    evil-N120-p98-chv12x10 & 120 & 545 & 20 & 24.526 & 24.526 & (1) & \cellcolor{lightgray} 20.000 & (2) \\
    \hline
    evil-N120-p98-myc5x24 & 120 & 236 & 48 & 52.607 & 48.000 & (22) & 48.000 & (16)\\
    \hline
    evil-N121-p98-myc11x11 & 121 & 508 & 22 & 26.397 & 26.397 & (1) & \cellcolor{lightgray} 22.000 & (1) \\
    \hline
    evil-N125-p98-s3m25x5 & 125 & 873 & 20 & 25.000 & 22.361 & (1) & 22.361 & (5) \\
    \hline
    evil-N138-p98-myc23x6 & 138 & 1242 & 12 & 15.177 & 15.177 & (1) & 15.177 & (3) \\
    \hline
    evil-N150-p98-myc5x30 & 150 & 338 & 60 & 65.121 & 60.000 & (32) & 60.000 & (43) \\
    \hline
    evil-N150-p98-s3m25x6 & 150 & 1102 & 24 & 30.000 & 26.833 & (2) & 26.833 & (8) \\
    \hline
    evil-N154-p98-myc11x14 & 154 & 701 & 28 & 33.596 & 33.596 & (1) & \cellcolor{lightgray} 28.000 & (2) \\
    \hline
    evil-N180-p98-chv12x15 & 180 & 944 & 30 & 36.788& 36.788 & (1) & \cellcolor{lightgray} 30.000 &(5)\\
    \hline
    evil-N184-p98-myc23x8 & 184 & 1764 & 16 & 20.235 & 20.235 & (2) & 20.235 & (7) \\
    \hline
	\end{tabular}
    %\end{adjustbox}
	\end{center}
	\label{table_4}
    \end{table}

    From the results for spin graphs presented in Table~\ref{table_4}, we observe that the values of BOUND 1 and BOUND 2 coincide with the stability numbers of considered graphs. This is particularly noteworthy as it indicates that we found the optimal stability number for the spin7 graph. Previously, only lower and upper bounds for~$\alpha(G)$ were reported in~\cite{Gaar2022}. These results highlight the potential for further research into SDP-based upper bounds for spin graphs, especially through the enhancement of the SDP by incorporating additional valid inequalities to refine the computation of the Lovász theta function~\eqref{theta_stable_set}.

    For the DIMACS instances, as shown in Table~\ref{table_4}, BOUND 1 already significantly improves upon the Lovász theta function as an upper bound on~$\alpha(G)$ for most of the considered graphs. BOUND 2 further enhances these bounds, offering notable improvements over BOUND 1 for several instances. In particular, BOUND 2 provides better integer bounds on~$\alpha(G)$ for three of the graphs. Additionally, the computations for both BOUND 1 and BOUND 2 were executed efficiently, with the only exception being the graph MANN\_a27.

    Furthermore, we examine the results for evil instances presented in Table~\ref{table_4}. As this is the first study of SDP-based bounds for these instances, we also review the values of the Lovász theta function~$\vartheta(G)$. Notably,~$\vartheta(G)$ consistently overestimates~$\alpha(G)$ across all instances, with differences ranging from approximately~$3$ to~$6.8$, indicating its lack of tightness as an upper bound. BOUND 1 provided a tighter upper bound on~$\alpha(G)$ than~$\vartheta(G)$ in four out of ten instances, and in two cases, it exactly matched~$\alpha(G)$. BOUND 2 demonstrates significant improvements over BOUND 1. Notably, for six out of the ten instances considered, BOUND 2 exactly matches the stability number~$\alpha(G)$. Interestingly, in instances where BOUND 2 does not exactly match~$\alpha(G)$, it provides the same value as BOUND 1. In particular, we observe that for instances evil-N138-p98-myc23x6 and evil-N184-p98-myc23x8, both BOUND 1 and BOUND 2 coincide with~$\vartheta(G)$. Therefore, it is a future research to explore whether incorporating some other inequalities into the SDP to compute the Lovász theta function might improve the bounds on stability numbers for these instances.

    Finally, we note that the computational effort required to calculate BOUND 1 and BOUND 2 was relatively modest. In particular, we obtained stability numbers for several instances in remarkably short running times for most instances. This not only highlights the quality of the results but also emphasizes the rapid computation of upper bounds on the stability numbers of evil graphs.

    \subsubsection{Third Set of Experiments}\label{third_set_stable_set}

    So far, we have presented computational results for relatively sparse instances. To ensure a more comprehensive analysis, we now perform a computational study including some denser instances. Specifically, we examine random graphs from the Erdős–Rényi model~$G(n, p)$ with~$n = 100$ and~$p \in \{0.10, 0.15, \ldots, 0.45, 0.50\}$.

    The main goal of this study is to evaluate the quality of the bounds BOUND 1 and BOUND 2 as the density of the graphs increases. Computational results, including the values of these bounds and the time required to compute them, are presented in Table~\ref{table_6}. The table follows the same structure as those from the second set of experiments, with additional information about the value of~$p$. Again, if the value of BOUND 2 surpasses BOUND 1 by at least one integer value, we highlight this result in gray. 

    \begin{table}[h!]
	\caption{Stable set results for random graphs on~$100$ vertices}
	\footnotesize
	\begin{center}
        %\begin{adjustbox}{width=\textwidth}
	\begin{tabular}{|l r r r c r | r r| r r|}
	\hline
    Graph &~$n$ &~$p$ &~$m$ &~$\alpha(G)$ &~$\vartheta(G)$ &
    BOUND 1&(time) & BOUND 2& (time)\\
    \hline
    rand\_100\_0.10 & 100 & 0.10 & 500 & 31 & 33.221 & 32.156 &  (5) & \cellcolor{lightgray} 31.279 & (18) \\
    \hline
    rand\_100\_0.15 & 100 & 0.15 & 757 & 24 & 26.276 & 25.740 & (3) & \cellcolor{lightgray} 24.598 & (50) \\
    \hline
    rand\_100\_0.20 & 100 & 0.20 & 1023 & 20 & 21.736 & 21.525 & (2) & \cellcolor{lightgray} 20.571 & (143) \\
    \hline
    rand\_100\_0.25 & 100 & 0.25 & 1265 & 17 & 19.078 & 18.843 & (2) & 18.074 & (225) \\
    \hline
    rand\_100\_0.30 & 100 & 0.30 & 1502 & 14 & 16.776 & 16.630 & (2) & \cellcolor{lightgray} 15.911 & (401) \\
    \hline
    rand\_100\_0.35 & 100 & 0.35 & 1790 & 13 & 14.413 & 14.274 & (3) & \cellcolor{lightgray} 13.750 & (655) \\
    \hline
    rand\_100\_0.40 & 100 & 0.40 & 1940 & 12 & 13.459 & 13.369 & (3) & \cellcolor{lightgray} 12.872 & (839) \\
    \hline
    rand\_100\_0.45 & 100 & 0.45 & 2224 & 10 & 11.705 & 11.609 & (3) & 11.343 & (886) \\
    \hline
    rand\_100\_0.50 & 100 & 0.50 & 2449 &  9 & 10.531 & 10.531 & (6) & 10.292 & (2317) \\
    \hline
	\end{tabular}
    %\end{adjustbox}
	\end{center}
	\label{table_6}
    \end{table}

Results from Table~\ref{table_6} show that BOUND 2 significantly improves upon BOUND 1 for sparse instances. Additionally, BOUND 2 yielded better integral solutions compared to BOUND 1 for instances with densities from~$0.30$ to~$0.40$. However, the computational time needed to obtain these bounds was quite high, as enumerating all cliques with at most~$5$-vertices and~$5$-cycles in denser graphs involves considerable computational effort.

\subsubsection{Analysis of Computational Demand}

We now provide details on the computational effort required to compute BOUND 1 and BOUND 2 for all instances analyzed in this computational study. While the reported times in Tables~\ref{table_1}, \ref{table_4}, and~\ref{table_6} offer an initial insight into the computational demands, we further examine the number of inequalities added to the SDP to compute the Lovász theta function for obtaining these bounds, as well as the number of iterations performed. These details are presented in Table~\ref{table_summary_1_1}. 

    \begin{table}[h!]
	\caption{Number of inequalities added and iterations performed for computing the bounds presented in Tables~\ref{table_1}, \ref{table_4}, and~\ref{table_6}}
	\footnotesize
    %\subfootnotesize
	\begin{center}
        %\begin{adjustbox}{width=\textwidth}
		\begin{tabular}{|l | r r | r r |}
		\hline
         & \multicolumn{2}{c|}{BOUND 1}  & \multicolumn{2}{c|}{BOUND 2} \\
         \cline{2-5}
		Graph & \# inequalities & \# iterations & \# inequalities & \# iterations \\
		\hline
        reg\_n100\_r4 & 2130 & 6 & 973 & 3 \\
        \hline
        reg\_n100\_r6 & 1263 & 4 & 1629 & 4 \\
        \hline
        reg\_n100\_r8 & 1004 & 3 & 1623 & 4 \\
        \hline
        reg\_n200\_r4 & 4642 & 7 & 2613 & 4 \\
        \hline
        reg\_n200\_r6 & 2692 & 4 & 3890 & 4 \\
        \hline
        reg\_n200\_r8 & 1608 & 3 & 3361 & 4 \\
        \hline
        reg\_n200\_r10 & 1285 & 2 & 3141 & 3 \\
        \hline
        rand\_n100\_p004 & 2474 & 6 & 126 & 1 \\
        \hline
        rand\_n100\_p006 & 1828 & 5 & 1027 & 2 \\
        \hline
        rand\_n100\_p008 & 1021 & 3 & 1517 & 3 \\
        \hline
        rand\_n100\_p010 & 1180 & 3 & 1226 & 3 \\
        \hline
        rand\_n200\_p002 & 7452 & 6 & 0 & 0 \\
        \hline
        rand\_n200\_p003 & 3739 & 4 & 1630 & 2 \\
        \hline
        rand\_n200\_p004 & 2384 & 3 & 3048 & 3 \\
        \hline
        rand\_n200\_p005 & 2323 & 3 & 2787 & 3 \\
        \hline
        torus\_5 & 100 & 1 & 50 & 1 \\
        \hline
        torus\_7 & 455 & 3 & 98 & 1 \\
        \hline
        torus\_9 & 2394 & 8 & 187 & 1 \\
        \hline
        torus\_11 & 4126 & 9 & 269 & 1 \\
        \hline
        torus\_13 & 7366 & 10 & 1282 & 2 \\
        \hline
        torus\_15 & 9901 & 10 & 2226 & 3 \\
        \hline
        spin5 & 1754 & 5 & 340 & 1 \\
        \hline
        spin7 & 10763 & 9 & 2210 & 2 \\
        \hline
        MANN\_a9 & 456 & 3 & 0 & 0 \\
        \hline
        MANN\_a27 & 10189 & 10 & 3162 & 3 \\
        \hline
        C125.9 & 880 & 2 & 1735 & 3  \\ 
        \hline
        C250.9 &  1413 & 1 & 1254 & 1 \\
        \hline
        sanr200\_0\_9 & 1313 & 1 & 1579 & 2  \\
        \hline
        evil-N120-p98-chv12x10 & 0 & 0 & 441 & 1 \\
        \hline
        evil-N120-p98-myc5x24 & 2486 & 5 & 383 & 2 \\
        \hline
        evil-N121-p98-myc11x11 & 0 & 0 & 363 & 1 \\
        \hline
        evil-N125-p98-s3m25x5 & 134 & 1 & 0 & 0 \\
        \hline
        evil-N138-p98-myc23x6 & 0 & 0 & 0 & 0\\
        \hline
        evil-N150-p98-myc5x30 & 3616 & 5 & 538 & 2 \\
        \hline
        evil-N150-p98-s3m25x6 & 210 & 1 & 0 & 0 \\
        \hline
        evil-N154-p98-myc11x14 & 0 & 0 & 465 & 1 \\
        \hline
        evil-N180-p98-chv12x15 & 0 & 0 & 665 & 1 \\
        \hline
        evil-N184-p98-myc23x8 & 0 & 0 & 0 & 0 \\
        \hline
        rand\_100\_0.10 &  1092 & 3 & 1339 & 3 \\
        \hline
        rand\_100\_0.15 & 633 & 2 & 1365 & 3 \\
        \hline
        rand\_100\_0.20 & 369 & 1 & 1146 & 3  \\
        \hline
        rand\_100\_0.25 & 380 & 1 & 884 & 2 \\
        \hline
        rand\_100\_0.30 & 253 & 1 & 713 & 2 \\
        \hline
        rand\_100\_0.35 & 268 & 1 & 532 & 2 \\
        \hline
        rand\_100\_0.40 & 214 & 1 & 533 & 2 \\
        \hline
        rand\_100\_0.45 &  223 & 1 & 270 & 1 \\
        \hline
        rand\_100\_0.50 & 180 & 1 & 169 & 1 \\
        \hline
		\end{tabular}
        %\end{adjustbox}
	\end{center}
	\label{table_summary_1_1}
    \end{table}
 
    From Table~\ref{table_summary_1_1}, we observe that the computation of BOUND 1 required more inequalities than BOUND 2 for $22$ instances, whereas BOUND 2 involved more inequalities for $23$ instances. At first glance, there seems to be no consistent pattern regarding which bound generally demands more inequalities. Upon closer examination, however, we find that even when the computation of BOUND 2 required more additional inequalities than BOUND 1, the total remained under~$4000$ across all instances, with at most~$4$ iterations needed. In contrast, the computation of BOUND 1 for~$7$ instances required addition of more that~$4000$ inequalities and, for some cases, required up to~$10$ iterations. This observation is interesting because, despite involving fewer inequalities and iterations, BOUND 2 typically produces stronger bounds than BOUND 1.

    To gain further insight into the inequalities added during computation, we provide a detailed breakdown of the specific inequalities incorporated into the SDP for computing the Lovász theta function. This analysis, aimed at identifying which inequalities were most frequently added, is presented in Table~\ref{table_summary_1_2}.

    \begin{table}[h!]
	\caption{Breakdown of added inequalities by type for computing the bounds presented in Tables~\ref{table_1}, \ref{table_4}, and~\ref{table_6}}
	\footnotesize
    \begin{center}
        %\begin{adjustbox}{width=\textwidth}
		\begin{tabular}{| l | r r r | r r r r r |}
		\hline
         & \multicolumn{3}{c|}{BOUND 1} & \multicolumn{5}{c|}{BOUND 2} \\
         \cline{2-9}
		Graph & \eqref{non-n} & \eqref{triangle_inequalities_1} & \eqref{triangle_inequalities_2} & \eqref{cliques_3} & \eqref{cliques_final} & \eqref{cycles_new} & \eqref{stable_set_2} & \eqref{stable_set_3}  \\
		\hline
        reg\_n100\_r4 & 250 & 1169 & 711 &  189 & 499 & 0 & 165 & 120\\
        \hline
        reg\_n100\_r6 & 252 & 788 & 223 &  203 & 581 & 9 & 600 & 236 \\
        \hline
        reg\_n100\_r8 & 297 & 579 & 128 & 368 & 538 & 39 & 624 & 54 \\
        \hline
        reg\_n200\_r4 & 458 & 2665 & 1519 & 418 & 1531 & 0 & 362 & 302 \\
        \hline
        reg\_n200\_r6 & 566 & 1600 & 526 & 456 & 1442 & 4 & 1362 & 626 \\
        \hline
        reg\_n200\_r8 & 512 & 1038 & 58 & 604 & 1163 & 59 & 1364 & 171 \\
        \hline
        reg\_n200\_r10 & 479 & 800 & 6 & 695 & 1059 & 151 & 1200 & 36 \\
        \hline
        rand\_n100\_p004 & 768 & 1019 & 687 & 13 & 30 & 1 & 43 & 39 \\
        \hline
        rand\_n100\_p006 & 409 & 898 & 521 & 115 & 351 & 19 & 327 & 215 \\
        \hline
        rand\_n100\_p008 & 249 & 567 & 205 & 339 & 516 & 66 & 525 & 71 \\
        \hline
        rand\_n100\_p010 & 447 & 516 & 217 & 204 & 403 & 65 & 456 & 98 \\
        \hline
        rand\_n200\_p002 & 3138 & 2331 & 1983 & 0 & 0 & 0 & 0 & 0 \\
        \hline
        rand\_n200\_p003 & 1206 & 1578 & 955 & 189 & 516 & 12 & 483 & 430 \\
        \hline
        rand\_n200\_p004 & 740 & 1200 & 444 & 554 & 1002 & 59 & 1021 & 412 \\
        \hline
        rand\_n200\_p005 & 989 & 995 & 339 & 562 & 937 & 104 & 1056 & 128 \\
        \hline
        torus\_5 & 0 & 50 & 50 & 0 & 50 & 0 & 0 & 0 \\
        \hline
        torus\_7 & 0 & 259 & 196 & 0 & 98 & 0 & 0 & 0 \\
        \hline
        torus\_9 & 4 & 1296 & 1094 & 25 & 162 & 0 & 0 & 0 \\
        \hline
        torus\_11 & 3 & 2152 & 1971 & 27 & 242 & 0 & 0 & 0 \\
        \hline
        torus\_13 & 1 & 3718 & 3647 & 606 & 676 & 0 & 0 & 0 \\
        \hline
        torus\_15 & 1 & 4950 & 4950 & 1279 & 947 & 0 & 0 & 0 \\
        \hline
        spin5 & 82 & 1152 & 520 & 93 & 250 & 0 & 0 & 0 \\
        \hline
        spin7 & 172 & 6091 & 4500 & 838 & 1372 & 0 & 0 & 0 \\
        \hline
        MANN\_a9 & 6 & 259 & 191 & 0 & 0 & 0 & 0 & 0 \\
        \hline
        MANN\_a27 &  502 & 8316 & 1371 & 1916 & 1246 & 0 & 0 & 0\\
        \hline
        C125.9 & 439 & 407 & 34 & 417 & 597 & 82 & 623 & 16  \\ 
        \hline
        C250.9 & 1296 & 117 & 0 & 397 & 500 & 3 & 354 & 0  \\
        \hline
        sanr200\_0\_9 & 912 & 400 & 1 & 483 & 609 & 17 & 470 & 0 \\
        \hline
        evil-N120-p98-chv12x10 & 0 & 0 & 0 & 201 & 240 & 0 & 0 & 0 \\
        \hline
        evil-N120-p98-myc5x24 & 644 & 1066 & 776 & 98 & 216 & 2 & 54 & 13 \\
        \hline
        evil-N121-p98-myc11x11 & 0 & 0 & 0 & 121 & 242 & 0 & 0 & 0 \\
        \hline
        evil-N125-p98-s3m25x5 & 9 & 125 & 0 & 0 &0 & 0 & 0 & 0 \\
        \hline
        evil-N138-p98-myc23x6 & 0 & 0 & 0 & 0 & 0 & 0 & 0 & 0 \\
        \hline
        evil-N150-p98-myc5x30 & 1021 & 1369 & 1226 & 146  & 318 & 5 & 58 & 11 \\
        \hline
        evil-N150-p98-s3m25x6 & 60 & 150 & 0 & 0 & 0 & 0 & 0 & 0 \\
        \hline
        evil-N154-p98-myc11x14 & 0 & 0 & 0 & 157 & 308 & 0 & 0 & 0\\
        \hline
        evil-N180-p98-chv12x15 & 0 & 0 & 0 & 305 & 360 & 0 & 0 & 0 \\
        \hline
        evil-N184-p98-myc23x8 & 0 & 0 & 0 & 0 & 0 & 0 & 0 & 0 \\
        \hline
        rand\_100\_0.10 & 437 & 492 & 163 & 272 & 450 & 75 & 499 & 43  \\
        \hline
        rand\_100\_0.15 & 377 & 253 & 3 & 373 & 528 & 34 & 430 & 0 \\
        \hline
        rand\_100\_0.20 & 277 & 92 & 0 & 361 & 479 & 3 & 303 & 0  \\
        \hline
        rand\_100\_0.25 & 297 & 83 & 0 & 320 & 398 & 0 & 166 & 0\\
        \hline
        rand\_100\_0.30 & 240 & 13 & 0 & 324 & 315 & 0 & 74 & 0  \\
        \hline
        rand\_100\_0.35 & 267 & 1 & 0 & 280 & 242 & 0 & 10 & 0 \\
        \hline
        rand\_100\_0.40 & 214 & 0 & 0 & 302 & 214 & 0 & 17 & 0 \\
        \hline
        rand\_100\_0.45 & 223 & 0 & 0 & 200 & 68 & 0 & 2 & 0  \\
        \hline
        rand\_100\_0.50 & 180 & 0 & 0 & 167 & 1 & 0 & 1 & 0 \\
        \hline
		\end{tabular}
        %\end{adjustbox}
	\end{center}
	\label{table_summary_1_2}
    \end{table}

 From Table~\ref{table_summary_1_2}, we observe that the computation of BOUND 1 generally involves more triangle inequalities~\eqref{triangle_inequalities_1} and~\eqref{triangle_inequalities_2} than non-negativity constraints~\eqref{non-n}, particularly for sparser instances. In denser instances, the addition of both triangle inequalities and non-negativity constraints is relatively modest. This is particularly the case for denser instances considered in the third set of experiments.

 Regarding the computation of BOUND 2, we note that the number of added inequalities~\eqref{cliques_final} valid for~$5$-cycles was very low, with at most~$151$ inequalities added for the instance reg\_n200\_r10. For inequalities~\eqref{cliques_3}, \eqref{cliques_final}, \eqref{stable_set_2}, and~\eqref{stable_set_3}, no clear pattern emerges. However, most instances required adding some inequalities valid for the join of cliques~\eqref{cliques_3} and~\eqref{cliques_final}, while inequalities valid for~$5$-cycles and one additional vertex~\eqref{stable_set_2} and~\eqref{stable_set_3} were not added in several instances. This indicates that no violations of these inequalities were found for those instances.

\subsubsection{Excluding Triangle Inequalities: Further Computational Study}\label{computational_result_excl_ti}

In the previous sections, we observed that once BOUND 1 is computed by incorporating non-negativity constraints and triangle inequalities into the SDP for computing the Lovász theta function, computing BOUND 2 becomes both efficient and straightforward. In particular, BOUND 2 shows a considerable improvement over BOUND 1 and provides strong upper bounds on~$\alpha(G)$ for most of the considered instances. However, our analysis of computational demand revealed an interesting pattern: the majority of the added inequalities required for computing BOUND 1, and thus BOUND 2, were triangle inequalities.

Motivated by this observation, we now conduct an additional study in which we exclude triangle inequalities from the computation of BOUND 1. To this end, we proceed as follows. After computing the optimal solution of the Lovász theta function~\eqref{theta_stable_set}, we examine the obtained solution and search for violations of the non-negativity constraints~\eqref{non-n}. We denote this bound as BOUND 1*. Then, we add  inequalities~\eqref{cliques_3},~\eqref{cliques_final}, \eqref{cycles_new}, \eqref{stable_set_2} and \eqref{stable_set_3} in the same manner as we did in our previous experiments. The resulting upper bound on~$\alpha(G)$ is denoted by BOUND 2*.

The primary objective of this experiment is to analyze the difference between the values of BOUND 2 and BOUND 2* and to assess the relative contribution of triangle inequalities to the overall quality of the computed bounds.

We consider instances from the previous experiments for which the computation of BOUND 1 required adding triangle inequalities and present the corresponding computational results in Table~\ref{table_excl_triangle_in}. In this table, we report for each instance the value of the Lovász theta function~$\vartheta(G)$, along with the already presented values of BOUND 1 and BOUND 2. Additionally, we include the newly computed bounds BOUND 1* and BOUND 2*, along with the computation times required for each bound.

\begin{table}[h!]
	\caption{Excluding triangle inequalities: Stable set results for selected instances considered in previous experiments}
	\footnotesize
	\begin{center}
        %\begin{adjustbox}{width=\textwidth}
		\begin{tabular}{|l r | rr| rr|rr|}
		\hline
    Graph &~$\vartheta(G)$  &  BOUND 1& BOUND 2 & BOUND 1* & (time) & BOUND 2* & (time) \\
		\hline
    reg\_n100\_r4  & 43.449  & 41.246  & 40.753 & 43.369 & (1) & 40.768 & (21) \\
        \hline
    reg\_n100\_r6 & 37.815 & 36.224 & 35.046 & 37.647 & (1) & 35.067 & (22) \\
    \hline
    reg\_n100\_r8 & 34.480 & 33.337 & 32.031 & 34.242 & (1) & 32.052 & (22) \\  
        \hline
    reg\_n200\_r4 & 87.759 & 83.498 &  82.292 & 87.646 & (1) & 82.279 & (103) \\ 
        \hline
    reg\_n200\_r6 & 79.276 & 76.047 &  73.738 & 79.027 & (2) & 73.718 & (120) \\
        \hline
    reg\_n200\_r8 & 70.790 & 69.110  & 66.743 & 70.566 & (2) & 66.766 & (73)  \\
        \hline
    reg\_n200\_r10& 66.418 & 65.142 & 62.714 & 66.168 & (2) & 62.653 & (77) \\
        \hline
    rand\_n100\_p004 & 46.067 & 45.032 & 45.006 & 45.955 & (1) & 45.012 & (12) \\
        \hline
    rand\_n100\_p006 & 40.361 & 38.909 & 38.435 & 40.182 & (1) & 38.452 & (14) \\
        \hline
    rand\_n100\_p008 & 34.847 & 33.575 & 32.431 & 34.690 & (1) & 32.428 & (27) \\
        \hline
    rand\_n100\_p010 & 34.020 & 32.934 & 32.126 & 33.790 & (1) & 32.151 & (20) \\
        \hline
    rand\_n200\_p002 & 95.778 & 95.044 & 95.044 & 95.962 & (8) & 95.003 & (54) \\
    \hline
    rand\_n200\_p003 & 83.662 & 81.560 &  81.064 & 83.372 & (2) & 81.061 & (101) \\
        \hline
    rand\_n200\_p004 & 73.908 & 71.654 & 69.854 & 73.615 & (2) & 69.825 & (122) \\
        \hline
    rand\_n200\_p005 & 69.039 & 67.313 & 65.509 & 68.664 & (3) & 65.487 & (133) \\
        \hline
    torus\_5 & 11.180 & 10.000 & 10.000 & 11.180 & (1) & 10.000 & (1)  \\
        \hline
    torus\_7 & 23.224 & 21.000 & 21.000 & 23.224 & (1) & 21.005 & (3) \\
        \hline
    torus\_9 & 39.241 & 36.000 & 36.000 & 39.241 & (1) &36.046 & (12) \\
        \hline
    torus\_11 & 59.249 & 55.022 & 55.019 & 59.249 & (1) & 55.051 & (26) \\
        \hline
    torus\_13 & 83.254 & 78.379 & 78.048 & 83.254 & (1) & 78.159 & (45)  \\
        \hline
    torus\_15 & 111.257 & 108.208 & 105.214 & 111.257 & (1) & 105.279 & (122) \\
        \hline
    spin5 & 55.902 & 50.000  & 50.000 & 55.902 & (1) & 55.002 & (18) \\
    \hline
    spin7 &162.566& 147.000 & 147.000 & 162.566 & (2) & 147.089 & (373) \\
    \hline
    MANN\_a9 & 17.475 & 17.220 & 17.220 & 17.475 & (1) & 17.090 & (1) \\
    \hline
    MANN\_a27 & 132.763 & 131.709 & 131.112 & 132.763 & (2) & 131.112 & (11) \\
    \hline
    C125.9 & 37.805 & 36.920 & 35.542 & 37.555 & (1) & 35.569 & (31) \\ 
    \hline
    C250.9 &  56.241 & 55.771 & 54.899 & 55.828 & (13) & 54.880 & (851) \\
    \hline
    sanr200\_0\_9 & 49.274 & 48.723 & 47.474 & 48.917 & (6) & 47.476 & (330) \\
    \hline
    evil-N120-p98-myc5x24 & 52.607 & 48.000 & 48.000 & 52.510 & (1) & 48.000 & (30) \\
    \hline
    evil-N125-p98-s3m25x5 & 25.000 & 22.361 & 22.361 & 25.000 & (1) & 22.532 & (138) \\
    \hline
    evil-N150-p98-myc5x30 & 65.121 & 60.000 & 60.000 & 64.983 & (1) & 60.000 & (100) \\
    \hline
    evil-N150-p98-s3m25x6 & 30.000 & 26.833 & 26.833 & 30.000 & (1) & 26.970 & (310) \\
    \hline
    rand\_100\_0.10 & 33.221 & 32.156 &  31.279 &  32.992 & (1) & 31.305 & (17)\\
    \hline
    rand\_100\_0.15 & 26.276 & 25.740 & 24.598 & 26.021 & (1) & 24.607 & (24)  \\
    \hline
    rand\_100\_0.20 & 21.736 & 21.525  & 20.571  & 21.601 & (1) & 20.572 & (127)\\
    \hline
    rand\_100\_0.25 & 19.078 & 18.843  & 18.074 &  18.882 & (2) & 18.080 & (192) \\
    \hline
    rand\_100\_0.30 & 16.776 & 16.630 & 15.911 &  16.641 & (2) & 15.920 & (338) \\
    \hline
    rand\_100\_0.35 & 14.413 & 14.274 &  13.750 & 14.274 & (3) & 13.750 & (550) \\
    \hline
		\end{tabular}
        %\end{adjustbox}
	\end{center}
	\label{table_excl_triangle_in}
    \end{table}

Results presented in Table~\ref{table_excl_triangle_in} reveal that adding the newly proposed inequalities \eqref{cliques_3}, \eqref{cliques_final}, \eqref{cycles_new}, \eqref{stable_set_2}, and~\eqref{stable_set_3}, along with the non-negativity constraints~\eqref{non-n}, into the SDP for computing the Lovász theta function yields strong upper bounds on the stability numbers of graphs. In particular, for all considered instances, BOUND 2* produces the same integer bound for~$\alpha(G)$ as BOUND 2. Furthermore, the computational times demonstrate that BOUND 2* was computed faster than BOUND 2 for almost all instances, highlighting the efficiency of this approach.

Overall, we conclude that although adding triangle inequalities can be done rather efficiently, their contribution to computing strong upper bounds on stability numbers of graphs becomes unnecessary when the inequalities proposed in this work, as well as non-negativity constraints, are included. 

We finalize this computational study by analyzing the computational effort required to compute the bounds BOUND 1* and BOUND 2*. To this end, we first present the number of inequalities added to SDP~\eqref{theta_stable_set} to achieve these bounds, along with the number of iterations performed. Subsequently, we provide a detailed breakdown of the types of inequalities added to SDP~\eqref{theta_stable_set} for the computation of these bounds. These findings are summarized in Tables~\ref{table_summary_2_1} and~\ref{table_summary_2_2}, respectively.

\begin{table}[h!]
	\caption{Number of inequalities added and iterations performed for computing the bounds presented in Table~\ref{table_excl_triangle_in}}
	\footnotesize

	\begin{center}
        %\begin{adjustbox}{width=\textwidth}
		\begin{tabular}{|l | r r | r r |}
		\hline
         & \multicolumn{2}{c|}{BOUND 1*}  & \multicolumn{2}{c|}{BOUND 2*} \\
         \cline{2-5}
		Graph & \# inequalities & \# iterations & \# inequalities & \# iterations \\
		\hline
        reg\_n100\_r4 & 174 & 1 & 2232 & 4 \\
        \hline
        reg\_n100\_r6 & 205 & 1 & 2333 & 4 \\
        \hline
        reg\_n100\_r8 & 251 & 1 & 2185 & 4 \\
        \hline
        reg\_n200\_r4 & 375 & 1 & 5297 & 5 \\
        \hline
        reg\_n200\_r6 & 511 & 1 & 5544 & 5 \\
        \hline
        reg\_n200\_r8 & 460 & 1 & 4363 & 4 \\
        \hline
        reg\_n200\_r10 & 438 & 1 & 4216 & 4 \\
        \hline
        rand\_n100\_p004 &  427 & 2 & 1471 & 3\\
        \hline
        rand\_n100\_p006 & 293 & 1 & 2185 & 3 \\
        \hline
        rand\_n100\_p008 & 193 & 1 & 2355 & 4 \\
        \hline
        rand\_n100\_p010 & 355 & 1 & 2075 & 3 \\
        \hline
        rand\_n200\_p002 & 2346 & 2 & 2090 & 3 \\
        \hline
        rand\_n200\_p003 & 845 & 1 & 4198 & 3 \\
        \hline
        rand\_n200\_p004 & 611 & 1 & 4653 & 4 \\
        \hline
        rand\_n200\_p005 & 827 & 1 & 4286 & 4  \\
        \hline
        torus\_5 &  0 & 0 & 200 & 1 \\
        \hline
        torus\_7 & 0 & 0 & 588 & 3 \\
        \hline
        torus\_9 & 0 & 0 & 1620 & 5 \\
        \hline
        torus\_11 & 0 & 0 & 2448 & 6 \\
        \hline
        torus\_13 & 0 & 0 & 3193 & 6 \\
        \hline
        torus\_15 & 0 & 0 & 4777 & 9 \\
        \hline
        spin5 & 0 & 0 & 2349 & 3 \\
        \hline
        spin7 & 0 & 0 & 6860 & 5 \\
        \hline
        MANN\_a9 & 0 & 0 & 72 & 1 \\
        \hline
        MANN\_a27 & 0 & 0 & 702 & 1 \\
        \hline
        C125.9 & 363 & 1 & 2365 & 3 \\ 
        \hline
        C250.9 & 1296 & 1 & 1468 & 1 \\
        \hline
        sanr200\_0\_9 & 912 & 1 & 1788 & 2  \\
        \hline
        evil-N120-p98-myc5x24 & 172 & 1 & 2793 & 4 \\
        \hline
        evil-N125-p98-s3m25x5 & 0 & 0 & 3511 & 7 \\
        \hline
        evil-N150-p98-myc5x30 & 284 & 1 & 4432 & 5 \\
        \hline
        evil-N150-p98-s3m25x6 & 0 & 0 & 5056 & 8 \\
        \hline
        rand\_100\_0.10 & 325 & 1 & 2095 & 3 \\
        \hline
        rand\_100\_0.15 & 328 & 1 & 1596 & 3 \\
        \hline
        rand\_100\_0.20 & 277 & 1 & 1193 & 3 \\
        \hline
        rand\_100\_0.25 & 297 & 1 & 926 & 2  \\
        \hline
        rand\_100\_0.30 & 240 & 1 & 759 & 2 \\
        \hline
        rand\_100\_0.35 & 267 & 1 & 532 & 2 \\
        \hline
		\end{tabular}
        %\end{adjustbox}
	\end{center}
	\label{table_summary_2_1}
    \end{table}

    From Table~\ref{table_summary_2_1}, we observe that incorporating the non-negativity constraints~\eqref{non-n} into the SDP to compute the Lovász theta function was highly efficient. Specifically, for computing BOUND 1, only a small number of non-negativity constraints were generally required, and the computations were completed in just a few iterations for nearly all instances.

    In contrast, computing BOUND 2* typically required adding more inequalities than were needed for computing BOUND 2, which consequently resulted in a greater number of iterations. However, the increase in the number of added inequalities was not substantial. As previously noted, the computations for BOUND 2* were generally completed in shorter times than BOUND 2.
    
    \begin{table}[h!]
	\caption{Breakdown of added inequalities by type for computing the bounds presented in Table~\ref{table_excl_triangle_in}}
	\footnotesize
	\begin{center}
        %\begin{adjustbox}{width=\textwidth}
		\begin{tabular}{| l | r | r r r r r |}
		\hline
         & \multicolumn{1}{c|}{BOUND 1*} & \multicolumn{5}{c|}{BOUND 2*} \\
         \cline{2-7}
		Graph & \eqref{non-n} & \eqref{cliques_3} & \eqref{cliques_final} & \eqref{cycles_new} & \eqref{stable_set_2} & \eqref{stable_set_3}  \\
		\hline
        reg\_n100\_r4 & 174 & 677 & 756 & 1 & 401 & 397 \\
        \hline
        reg\_n100\_r6 & 205 & 521 & 674 & 15 & 700 & 423\\
        \hline
        reg\_n100\_r8 & 251 & 504 & 691 & 80 & 694 & 216 \\
        \hline
        reg\_n200\_r4 & 375 & 1609 & 1988 & 0 & 875 & 825 \\
        \hline
        reg\_n200\_r6 & 511 & 1120 & 1781 & 5 & 1683 & 955 \\
        \hline
        reg\_n200\_r8 & 460 & 931 & 1375 & 73 & 1470 & 514 \\
        \hline
        reg\_n200\_r10 & 438 & 926 & 1295 & 231 & 1364 & 400 \\
        \hline
        rand\_n100\_p004 & 427 & 357 & 437 & 14 & 344 & 319 \\
        \hline
        rand\_n100\_p006 & 295 & 483 & 600 & 44 & 596 & 462 \\
        \hline
        rand\_n100\_p008 & 193 & 516 & 679 & 167 & 672 & 321 \\
        \hline
        rand\_n100\_p010 & 355 & 434 & 499 & 178 & 596 & 368 \\
        \hline
        rand\_n200\_p002 & 2346 & 550 & 614 & 3 & 466 & 457 \\
        \hline
        rand\_n200\_p003 & 845 & 891 & 1042 & 51 & 1174 & 1040 \\
        \hline
        rand\_n200\_p004 & 611 & 1024 & 1359 & 153 & 1360 & 757 \\
        \hline
        rand\_n200\_p005 & 827 & 924 & 1282 & 213 & 1305 & 562 \\
        \hline
        torus\_5 & 0 & 50 & 50 & 0 & 50 & 50 \\
        \hline
        torus\_7 & 0 & 294 & 294 & 0 & 0 & 0 \\
        \hline
        torus\_9 & 0 & 810 & 810 & 0 & 0 & 0 \\
        \hline
        torus\_11 & 0 & 1310 & 1138 & 0 & 0 & 0 \\
        \hline
        torus\_13 &0 & 1757 & 1436 & 0 & 0 & 0 \\
        \hline
        torus\_15 &  0 & 3371 & 1406 & 0 & 0 & 0 \\
        \hline
        spin5 & 0 & 548 & 748 & 0 & 521 & 532 \\
        \hline
        spin7 & 0 & 3430 & 3430 & 0 & 0 & 0 \\
        \hline
        MANN\_a9 & 0  & 72 & 0 & 0 & 0 & 0 \\
        \hline
        MANN\_a27 & 0 & 702 & 0 & 0 & 0 & 0 \\
        \hline
        C125.9 & 363 & 562 & 685 & 165 & 700 & 253 \\ 
        \hline
        C250.9 & 1296 & 464 & 500 & 4 & 500 & 0 \\
        \hline
        sanr200\_0\_9 & 912 & 514 & 690 & 37 & 547 & 0  \\
        \hline
        evil-N120-p98-myc5x24 & 172 & 787 & 850 & 15 & 591 & 550 \\
        \hline
        evil-N125-p98-s3m25x5 & 0 & 1740 & 1700 & 0 & 71 & 0 \\
        \hline
        evil-N150-p98-myc5x30 & 284 & 1221 & 1303 & 21 & 962 & 925 \\
        \hline
        evil-N150-p98-s3m25x6 & 0 & 2361 & 2386 & 0 & 309 & 0 \\
        \hline
        rand\_100\_0.10 & 325 & 452 & 542 & 200 & 600 & 301 \\
        \hline
        rand\_100\_0.15 & 328 & 434 & 559 & 93 & 496 & 14 \\
        \hline
        rand\_100\_0.20 & 277 & 373 & 484 & 7 & 329 & 0 \\
        \hline
        rand\_100\_0.25 & 297 & 312 & 400 & 0 & 214 & 0 \\
        \hline
        rand\_100\_0.30 &  240 & 324 & 323 & 0 & 112 & 0 \\
        \hline
        rand\_100\_0.35 & 267 & 280 & 242 & 0 & 10 & 0 \\
        \hline
		\end{tabular}
        %\end{adjustbox}
	\end{center}
	\label{table_summary_2_2}
    \end{table}

    Finally, from Table~\ref{table_summary_2_2}, we observe that for almost all instances, except for torus graphs, inequalities of types \eqref{cliques_3}, \eqref{cliques_final}, \eqref{cycles_new}, \eqref{stable_set_2}, and \eqref{stable_set_3} were added, with the smallest number being of type \eqref{cycles_new}. Altogether, this demonstrates that all these inequalities contributed to obtaining strong upper bounds on the stability numbers of graphs.

    \subsection{The Graph Coloring Problem}\label{results_coloring}

We now conduct a computational study to strengthen the Lovász theta function~$\vartheta(\overline{G})$ towards~$\chi(G)$. Our approach mirrors the methodology used for the stable set problem. Hence, we iteratively incorporate inequalities~\eqref{non-n-col}, \eqref{triangle_inequalities_coloring}, \eqref{clique_vertex_four}, \eqref{coloring_1}, and~\eqref{coloring_2} into the SDP for computing the Lovász theta function~\eqref{theta_coloring}, in a two-phase procedure:

\begin{itemize}
\item
\textbf{Phase 1:} Starting from the optimal solution of the formulation~\eqref{theta_coloring}, in each iteration, we add all violations of the non-negativity constraints~\eqref{non-n-col} and up to~$2n$ of the largest violations greater than~$0.025$ of the triangle inequalities~\eqref{triangle_inequalities_coloring}. This process continues until there are fewer than~$n$ violations or until the number of iterations reaches~$10$. The resulting lower bound on the chromatic number~$\chi(G)$ is denoted by BOUND 1.
\item
\textbf{Phase 2:} Using the optimal solution from the SDP to compute BOUND 1, we then address violations of inequalities valid for a clique and one additional vertex~\eqref{clique_vertex_four} as well as the proposed inequalities~\eqref{coloring_1} and \eqref{coloring_2}. For this purpose, we enumerate all cliques with at most~$5$ vertices and consider induced cycles of length~$5$. In each iteration, we add up to~$2n$ of the largest violations greater than~$0.025$ per inequality type. This phase continues until there are fewer than~$n$ violations or until the number of iterations reaches 10. We denote the resulting lower bound on~$\chi(G)$ by BOUND 2.
\end{itemize}

\subsubsection{First Set of Experiments}\label{first_experiments_coloring}

For the first set of experiments, we consider some of the instances analyzed in~\cite{Gaa:20}. Specifically, we examine a mug graph on~$88$ vertices, several Mycielski graphs, and several Full Insertion graphs, which are a generalization of Mycielski graphs. Mycielski graphs, introduced by Tomescu~\cite{Tomescu1968}, are particularly interesting because they are challenging to solve. Since they are triangle-free, their clique number is~$2$, but their chromatic number increases with problem size.
    
    The primary objective of this computational study is to evaluate the improvement of BOUND 2 over BOUND 1 as a lower bound on~$\chi(G)$ and to compare our results with those obtained in~\cite{Gaa:20}. Additionally, we aim to investigate the computational times required to compute these bounds.

    The results for this first set of experiments are given in Table~\ref{table_col_1}. The presentation of these results follows the same format as the data for the stable set problem. Thus, the first columns of Table~\ref{table_col_1} provide general information about the instances, i.e., the name of the instance, the number of vertices~$n$, the number of edges~$m$, the value of~$\chi(G)$, and the value of the Lovász theta function~$\vartheta(\overline{G})$. Then, we report the computed bounds BOUND 1 and BOUND 2, along with the bounds from~\cite{Gaa:20}, denoted as GR bounds. For each instance, the best computed bound on~$\chi(G)$ is highlighted in bold. Additionally, if BOUND 2 yields a better integer bound on~$\chi(G)$ than bound GR, the corresponding cell is shaded in gray. Finally, we present the computational times (time) required for the computation of bounds BOUND 1 and BOUND 2, as well as the times for the computation of the GR bound, as reported in \cite{Gaa:20}. Again, we note that our computational setup differs from that in~\cite{Gaa:20}. However, we include the reported computational times from~\cite{Gaa:20} to provide context for the computational effort involved.

    \begin{table}[h!]
	\caption{Graph coloring results for some instances considered in~\cite{Gaa:20}}
	\footnotesize
	\begin{center}
        %\begin{adjustbox}{width=\textwidth}
		\begin{tabular}{|l r r c c|  r r|rr|rr|}
		\hline
		Graph & $n$ & $m$ & $\chi(G)$ &
        $\vartheta(\overline{G})$ & BOUND 1 & (time)& BOUND 2 &(time) & GR &(time)  \\
		\hline
        myciel5 & 47 & 236 & 6 & 2.639 & 3.093 & (3) & 3.468 & (17) & \textbf{3.510} & (4240) \\
        \hline
        myciel6 & 95 & 755 & 7 & 2.734 & 3.253 & (21) & \textbf{3.622} & (406) & 3.534 & (1540) \\
        \hline
        mug88\_1 & 88 & 146 & 4 & 3.000 & 3.001 & (3) & 3.001 & (1) & \textbf{3.022} & (4709) \\
        \hline
        1\_FullIns\_4 & 93 & 593 & 5 & 3.124 & 3.487 & (3) & 3.837 & (23) & \textbf{3.939} & (7220) \\
        \hline
        2\_FullIns\_4 & 212 & 1621 & 6 & 4.056 & 4.343 & (4) & 4.670 & (17) & \textbf{4.700} & (10106) \\
        \hline        
		\end{tabular}
        %\end{adjustbox}
	\end{center}
	\label{table_col_1}
    \end{table}

    The results in Table~\ref{table_col_1} show that BOUND 1 consistently improves upon~$\vartheta(\overline{G})$ as a lower bound on~$\chi(G)$ across all five instances, with BOUND 2 offering additional enhancements over BOUND 1 in four out of five cases. When comparing BOUND 2 as a lower bound on~$\chi(G)$ to the GR bounds, BOUND 2 often comes very close to the GR results and for one instance, BOUND 2 actually surpasses GR. Nevertheless, for all considered instances, BOUND 2 provides the same integer bounds on~$\chi(G)$ as the GR bound.

    An important consideration is the computational time required for each method. While the GR bounds are often slightly superior, they come at a significantly higher computational cost. In contrast, for four out of five instances, BOUND 2 was computed in approximately~$20$ seconds. This efficiency highlights the advantage of BOUND 2, providing competitive bounds with much less computational effort.
    
\subsubsection{Second Set of Experiments}\label{second_experiments_coloring}

In our second set of experiments, we consider some instances from the literature:
    \begin{itemize}
    \item DSJC graphs: DSJC graphs are standard~$G(n, p)$ random graphs. They were initially used in the paper by Johnson, Aragon, McGeoch, and Schevon~\cite{Johnson1991} and were also included in the second DIMACS challenge~\cite{Johnson1996}. We consider two sparse DSJC graphs.
    \item Full Insertion graphs: We explore additional Full Insertion graphs that were not considered in~\cite{Gaa:20}. 
    \item Queen graphs: The~$n \times n$ queen graph is a graph on~$n^2$ vertices in which each vertex represents a square in an~$n \times n$ chessboard, and each edge corresponds to a legal move by a queen.
    \item Random graphs from~\cite{Dukanovic2007}: We examine some of the random graphs considered in~\cite{Dukanovic2007}, where~$\vartheta(\overline{G})$ was strengthened towards~$\chi(G)$ by adding the non-negativity constraints as well as triangle inequalities. Hence, bounds presented in~\cite{Dukanovic2007} should generally coincide with our BOUND 1, so we do not report their results separately. Moreover, we do not know the exact values of~$\chi(G)$ for these instances, hence these values are not reported.
    \end{itemize}

    In this computational study, we compare the values of bounds~$\vartheta(\overline{G})$, BOUND 1, and BOUND 2 as lower bounds on~$\chi(G)$ and assess the computational effort required for the computations of BOUND 1 and BOUND 2. The results are presented in Table~\ref{table_col_2}.

    \begin{table}[h!]
	\caption{Graph coloring results for selected instances from the literature}
	\footnotesize
	\begin{center}
        %\begin{adjustbox}{width=\textwidth}
		\begin{tabular}{|l r r c c|  r r|rr|}
		\hline
		Graph & $n$ &$m$ &$\chi(G)$ &
        $\vartheta(\overline{G})$ & BOUND 1 & (time)& BOUND 2 &(time) \\
		\hline
        dsjc125.1 & 125 & 736 & 5 & 4.106 & 4.218 & (2)& 4.430 & (14)\\ 
        \hline
        dsjc250.1 & 250 & 3218 & 8 & 4.906 & 4.939 &(12) & \cellcolor{lightgray} 5.040 & (457)\\
        \hline
        3\_FullIns\_3 & 80 & 346 & 6 & 5.016 & 5.194&(1)& 5.194&(1) \\
        \hline
        4\_FullIns\_3 & 114 & 541 & 7 & 6.010 &6.010 &(1)& 6.309&(1) \\
        \hline
        5\_FullIns\_3 & 154 & 792 & 8 & 7.007 &7.007&(1)& 7.267&(2) \\
        \hline
        Queen\_8\_8 & 64 & 728 & 9 & 8.000 & 8.000 &(1)& 8.000&(7) \\
        \hline
        Queen\_9\_9 & 81 & 1056 & 10 & 9.000 &9.000 &(1)& 9.000&(25) \\
        \hline
        Queen\_10\_10 & 100 & 1470 & 11 & 10.000 & 10.000&(1) & 10.000&(62) \\
        \hline
        G100\_25 & 100 & 1240 & - & 5.823 &5.867&(2)&\cellcolor{lightgray} 6.235&(76) \\
        \hline
        G150\_25 & 150 & 2802 & - & 6.864& 6.918&(6)&\cellcolor{lightgray} {7.184}&(529) \\
        \hline
        G200\_1 & 200 & 2047& - &4.447 & 4.473&(10) & 4.600&(178) \\
        \hline
        G250\_1 & 250 & 3149& - &4.805 & 4.831&(12) & 4.928&(512) \\
        \hline
		\end{tabular}
        %\end{adjustbox}
	\end{center}
	\label{table_col_2}
    \end{table}

    The results in Table~\ref{table_col_2} indicate that BOUND 1 improves upon~$\vartheta(\overline{G})$ as a lower bound on~$\chi(G)$ for seven out of twelve instances. BOUND 2, in comparison, provides better results for eight instances and surpasses BOUND 1 by one integer value in three cases, demonstrating its effectiveness in tightening lower bounds.

    However, BOUND 2 requires more computational effort, especially for instances where it improves upon BOUND 1 by one integer value. Interestingly, for all queen graphs considered, both bounds align with the value of the Lovász theta function~$\vartheta(\overline{G})$, suggesting that the added inequalities do not enhance this bound. Further research could explore whether other inequalities might improve~$\vartheta(\overline{G})$ as a lower bound on~$\chi(G)$.

\subsubsection{Third Set of Experiments}\label{third_experiments_coloring}

We now extend our computational study by considering random graphs, as previously investigated in the context of strengthening the Lovász theta function towards the stability numbers of graphs in Section~\ref{third_set_stable_set}. Hence, we examine random graphs from the Erdős–Rényi model~$G(n, p)$ with~$n = 100$ and~$p \in \{0.10, 0.15, \ldots, 0.45, 0.50\}$. 

The objective of this study is to evaluate the extent to which BOUND 2 provides tighter bounds on the chromatic number of these graphs compared to BOUND 1 and the Lovász theta function. The results are summarized in Table~\ref{table_col_3}. If BOUND 2 yields a better integral bound on~$\chi(G)$ than BOUND 1, we highlight the corresponding cell in gray.

\begin{table}[h!]
	\caption{Graph coloring results for random graphs on~$100$ vertices}
	\footnotesize
	\begin{center}
        %\begin{adjustbox}{width=\textwidth}
	\begin{tabular}{|l r r c r | r r| r r|}
	\hline
    Graph &~$n$ &~$p$ &~$m$ &~$\vartheta(\overline{G})$ &
    BOUND 1&(time) & BOUND 2& (time)\\
    \hline
    rand\_100\_0.10 & 100 & 0.10 & 500 & 4.000 & 4.000 & (1) & 4.000 & (1) \\
    \hline
    rand\_100\_0.15 & 100 & 0.15 & 757 & 4.373 & 4.413 & (2) & 4.708 & (47) \\
    \hline
    rand\_100\_0.20 & 100 & 0.20 & 1023 & 5.178 & 5.208 & (2) & 5.497 & (47) \\
    \hline
    rand\_100\_0.25 & 100 & 0.25 & 1265 & 5.859 & 5.908 & (2) & \cellcolor{lightgray} 6.287 & (77)  \\
    \hline
    rand\_100\_0.30 & 100 & 0.30 & 1502 & 6.466 & 6.503 & (2) & 6.915 & (104) \\
    \hline
    rand\_100\_0.35 & 100 & 0.35 & 1790 & 7.589 & 7.640 & (2) & \cellcolor{lightgray} 8.111 & (152) \\
    \hline
    rand\_100\_0.40 & 100 & 0.40 & 1940 & 8.030 & 8.088 & (2) & 8.629 & (216) \\
    \hline
    rand\_100\_0.45 & 100 & 0.45 & 2224 & 9.396 & 9.467 & (3) & \cellcolor{lightgray} 10.015 & (247) \\
    \hline
    rand\_100\_0.50 & 100 & 0.50 & 2449 & 10.321 & 10.399 & (6) & 10.900 & (716)  \\
    \hline
	\end{tabular}
    %\end{adjustbox}
	\end{center}
	\label{table_col_3}
    \end{table}

The results in Table~\ref{table_col_3} show that BOUND 2 improves upon BOUND 1 in eight out of nine instances. Additionally, for three instances, we obtained even tighter integer bounds on the chromatic numbers of graphs. Interestingly, these improvements were achieved for instances with densities of~$25\%$,~$35\%$, and~$45\%$. This is particularly significant because, as Table~\ref{table_col_3} indicates, adding non-negativity constraints and triangle inequalities provides only marginal improvement over the Lovász theta function for random graphs. In contrast, the proposed inequalities prove to be highly effective in strengthening lower bounds on chromatic numbers of graphs, even for denser random graphs.

\subsubsection{Analysis of Computational Demand}

We conclude our computational study by analyzing the effort required to compute bounds BOUND 1 and BOUND 2 presented in this section. For this purpose, we analyze only instances where BOUND 2 presents an improvement over the Lovász theta function. Table~\ref{table_summary_col_1} provides details on the number of iterations and inequalities added into the SDP for computing the Lovász theta function. Additionally, Table~\ref{table_summary_col_2} presents a breakdown of the added inequalities by type for computing the respective bounds.

  \begin{table}[h!]
	\caption{Number of inequalities added and iterations performed for computing the bounds presented in Tables~\ref{table_col_1}, \ref{table_col_2}, and~\ref{table_col_3}}
	\footnotesize
	\begin{center}
        %\begin{adjustbox}{width=\textwidth}
		\begin{tabular}{|l | r r | r r |}
		\hline
         & \multicolumn{2}{c|}{BOUND 1}  & \multicolumn{2}{c|}{BOUND 2} \\
         \cline{2-5}
		Graph & \# inequalities & \# iterations & \# inequalities & \# iterations \\
		\hline
        myciel5 &  452 & 5 & 1381 & 8 \\
        \hline
        myciel6 & 1330 & 7 & 4509 & 10 \\
        \hline
        mug88\_1 & 232 & 2 & 0 & 0 \\
        \hline
        1\_FullIns\_4 & 664 & 3 & 2121 & 5 \\
        \hline
        2\_FullIns\_4 & 424 & 1 & 1129 & 1 \\
        \hline    
        dsjc125.1 & 780 & 1 & 1076 & 3 \\ 
        \hline
        dsjc250.1 & 1144 & 1 & 2075 & 4 \\
        \hline
        3\_FullIns\_3 & 115 & 1 & 0 & 0 \\
        \hline
        4\_FullIns\_3 & 0 & 0 & 119 & 1 \\
        \hline
        5\_FullIns\_3 & 0 & 0 & 184 & 1 \\
        \hline
        G100\_25 & 311 & 1 & 1684 & 8  \\
        \hline
        G150\_25 & 536 & 1 & 1954 & 7 \\
        \hline
        G200\_1 & 815 & 1 & 2512 & 5 \\
        \hline
        G250\_1 & 1043 & 1 & 2076 & 4 \\
        \hline
        rand\_100\_0.10 & 0 & 0 & 308 & 2 \\
        \hline
        rand\_100\_0.15 & 570 & 2 & 2232 & 7 \\
        \hline
        rand\_100\_0.20 & 330 & 1 & 1708 & 7 \\
        \hline
        rand\_100\_0.25 & 312 & 1 & 1789 & 8 \\
        \hline
        rand\_100\_0.30 & 207 & 1 & 1609 & 8 \\
        \hline
        rand\_100\_0.35 & 206 & 1 & 1604 & 8 \\
        \hline
        rand\_100\_0.40 & 215 & 1 & 1773 & 9 \\
        \hline
        rand\_100\_0.45 & 212 & 1 & 1573 & 8 \\
        \hline
        rand\_100\_0.50 & 217 & 1 & 1400 & 7 \\
        \hline
		\end{tabular}
        %\end{adjustbox}
	\end{center}
	\label{table_summary_col_1}
    \end{table}

Table~\ref{table_summary_col_1} shows that for~$21$ instances where BOUND 2 improves upon BOUND 1, fewer iterations were required to obtain BOUND 1 than BOUND 2. In other words, adding non-negativity constraints and triangle inequalities was done highly efficiently, often requiring only a single iteration. Consequently, relatively few of these inequalities needed to be added, which aligns with the computational times reported in Tables~\ref{table_col_1}, \ref{table_col_2}, and~\ref{table_col_3}.

In contrast, incorporating the newly proposed inequalities generally required significantly more iterations. This indicates that a substantial number of violated constraints were identified and added into the SDP to compute the Lovász theta function. This observation is also reflected in the computational times presented in Tables~\ref{table_col_1}, \ref{table_col_2}, and~\ref{table_col_3}. Moreover, it aligns with the fact that BOUND 2 sometimes yields substantial improvements over BOUND 1.

\begin{table}[h!]
	\caption{Breakdown of added inequalities by type for computing the bounds presented in Tables~\ref{table_col_1}, \ref{table_col_2}, and~\ref{table_col_3}}
	\footnotesize
	\begin{center}
        %\begin{adjustbox}{width=\textwidth}
		\begin{tabular}{| l | r  r | r r r  |}
		\hline
         & \multicolumn{2}{c|}{BOUND 1} & \multicolumn{3}{c|}{BOUND 2} \\
         \cline{2-6}
		Graph & \eqref{non-n-col} & \eqref{triangle_inequalities_coloring} & \eqref{clique_vertex_four} & \eqref{coloring_1} & \eqref{coloring_2}   \\
		\hline
        myciel5 & 0 & 452 & 52 & 736 & 593   \\
        \hline
        myciel6 & 0 & 1330 & 329 & 2090 & 2090  \\
        \hline
        mug88\_1 & 30 & 202 & 0 & 0 & 0 \\
        \hline
        1\_FullIns\_4 & 106 & 558 & 598 & 713 & 810  \\
        \hline
        2\_FullIns\_4 & 0 & 424 & 349 & 356 & 424  \\
        \hline    
        dsjc125.1 & 530 & 250 & 682 & 70 & 324 \\ 
        \hline
        dsjc250.1 & 1107 & 37 & 2000 & 26 & 49 \\
        \hline
        3\_FullIns\_3 & 5 & 110 & 0 & 0 & 0\\
        \hline
        4\_FullIns\_3 & 0 & 0 & 42 & 77 & 0 \\
        \hline
        5\_FullIns\_3 & 0 & 0 & 56 & 128 & 0 \\
        \hline
        G100\_25 & 275 & 36 & 1513 & 27 & 144 \\
        \hline
        G150\_25 & 535 & 1 & 1952 & 0 & 2 \\
        \hline
        G200\_1 & 685 & 130 & 1970 & 114 & 428 \\
        \hline
        G250\_1 & 1049 & 24 & 1997 & 19 & 60 \\
        \hline
        rand\_100\_0.10 & 0 & 0 & 219 & 47 & 42  \\
        \hline
        rand\_100\_0.15 & 311 & 259 & 1367 & 262 & 603  \\
        \hline
        rand\_100\_0.20 & 254 & 76 & 1370 & 66 & 272  \\
        \hline
        rand\_100\_0.25 & 289 & 23 & 1518 & 50 & 221  \\
        \hline
        rand\_100\_0.30 & 203 & 4 & 1517 & 9 & 83  \\
        \hline
        rand\_100\_0.35 & 206 & 0 & 1600 & 1 & 3  \\
        \hline
        rand\_100\_0.40 & 215 & 0 & 1725 & 0 & 0  \\
        \hline
        rand\_100\_0.45 & 212 & 0 & 1573 & 0 & 0  \\
        \hline
        rand\_100\_0.50 & 217 & 0 & 1400 & 0 & 0  \\
        \hline
		\end{tabular}
        %\end{adjustbox}
	\end{center}
	\label{table_summary_col_2}
    \end{table}

When we analyze the breakdown of added inequalities by type presented in Table~\ref{table_summary_col_2}, we find that for the~$20$ instances where BOUND 1 improves upon the Lovász theta function, more non-negativity constraints than triangle inequalities were added in~$14$ cases, while the opposite was true in~$6$ cases. This aligns with prior observations from~\cite{Dukanovic2007}, which suggest that adding triangle inequalities does not significantly enhance the Lovász theta function.

From Table~\ref{table_summary_col_2}, we also observe that for all considered random graphs, a substantial number of newly proposed inequalities valid for cliques and one vertex, i.e., of type~\eqref{clique_vertex_four} were added. This highlights the importance of incorporating these inequalities into the SDP when computing the Lovász theta function to obtain strong bounds on the chromatic number of random graphs.

For several other instances, such as the Mycielski graph 6 and 1\_FullIns\_4, more inequalities of types~\eqref{coloring_1} and \eqref{coloring_2} which are valid for odd cycles were added. This is particularly interesting, as it demonstrates that to achieve tight bounds, all three proposed inequalities~\eqref{clique_vertex_four}, \eqref{coloring_1}, and \eqref{coloring_2} should be considered.

Furthermore, we observe that although triangle inequalities were added into the SDP for computing the Lovász theta function, violations of the newly proposed inequalities were always detected. This confirms our earlier observation that triangle inequalities do not imply the inequalities introduced in Section~\ref{inequalities_coloring}.

Finally, we note that in the previous section, we applied a similar two-phase approach to strengthen the Lovász theta function towards the stability number of a graph. In that case, a significant number of triangle inequalities were added for many instances, prompting a further computational study to assess the impact of excluding them. However, as shown in Tables~\ref{table_summary_col_1} and~\ref{table_summary_col_2}, this is not the case for the graph coloring problem. In almost all instances, only a modest number of triangle inequalities were added, making further computational analysis unnecessary.

 \subsection{Summary of Results}\label{summary}

    In Sections~\ref{results_stable_set} and~\ref{results_coloring}, we performed an extensive computational study of bounds for stability and chromatic numbers by incorporating valid inequalities into the SDP for the computation of the Lovász theta function. We assessed these bounds across various instances from the literature and provided a detailed analysis.

    Overall, our findings indicate that incorporating non-negativity constraints and triangle inequalities into the semidefinite program for the Lovász theta function enables an efficient transition to computing BOUND 2. In addition to being easy and fast to compute, BOUND 2 yields significantly stronger results than BOUND 1 for several classes of instances, often producing better integer solutions, as highlighted in Section~\ref{results_stable_set}. Interestingly, for the stable set problem, strong bounds on $\alpha(G)$ can also be achieved by computing BOUND 1* and BOUND 2*, which omit the triangle inequalities.
    
    Furthermore, it is worth mentioning that we limited our iterative computational approach to a maximum of ten iterations, but only a few iterations were required for most instances, further demonstrating the practical value of this approach for deriving upper bounds on $\alpha(G)$ and lower bounds on $\chi(G)$.

    We observe that the effectiveness of the computed bounds on $\alpha(G)$ and $\chi(G)$ generally depends on the size and sparsity of the instance. The best results are obtained for smaller and sparser instances, while performance tends to decrease for larger and denser ones. Interestingly, for certain evil instances in Section~\ref{results_stable_set} as well as queen graphs in Section~\ref{results_coloring}, adding inequalities into the SDP did not enhance the Lovász theta function. Despite this, our bounds remain highly competitive with those found in the literature and can be computed in relatively short running times.
	
	\section{Conclusion}\label{Conclusion}

    Several questions remain open for future investigation. First, it would be valuable to explore which additional inequalities could be added into the SDP for the computation of the Lovász theta function to improve upper bounds on $\alpha(G)$ for certain evil graphs and lower bounds on $\chi(G)$ for queen graphs. Additionally, our research thus far has concentrated on subgraphs containing cliques, odd cycles, and odd antiholes. It is an open question to identify other easily recognizable structures and examine whether incorporating valid inequalities for these structures could further enhance the Lovász theta function with respect to $\alpha(G)$ and $\chi(G)$. Regarding the new inequalities introduced in this work, it would be beneficial to determine which additional inequalities are needed to fully characterize the underlying polytopes, as we have done for subgraphs containing cliques in~\cite{PuchRen2023}. Finally, given the strength of the proposed bounds,  developing an exact algorithm for the stable set problem within a branch-and-bound framework that utilizes these bounds would be of significant interest.
	
    \section*{Disclosure statement}

    The authors report there are no competing interests to declare.

    \section*{Data availability statement}

    The program code associated with this paper is available as ancillary files from the arXiv page of this paper. Additionally, the source code and data are available upon request from the authors.

    \section*{Acknowledgments}

    We gratefully acknowledge Adam Letchford for his insightful comments and helpful suggestions.

    \section*{Funding} This research was funded in whole by the Austrian Science Fund (FWF) [10.55776/DOC78]. For the purposes of open access, the authors have applied a CC BY public copyright license to all author-accepted manuscript versions resulting from this submission.
	
	\afterpage{\clearpage}
	
	\bibliographystyle{plain}
	\bibliography{References}
	
\end{document}